# BOUNDARY SINGULARITY OF POISSON AND HARMONIC BERGMAN KERNELS

Miroslav Engliš

ABSTRACT. We give a complete description of the boundary behaviour of the Poisson kernel and the harmonic Bergman kernel of a bounded domain with smooth boundary, which in some sense is an analogue of the similar description for the usual Bergman kernel on a strictly pseudoconvex domain due to Fefferman. Our main tool is the Boutet de Monvel calculus of pseudodifferential boundary operators, and in fact we describe the boundary singularity of a general potential, trace or singular Green operator from that calculus.

## 1. INTRODUCTION

Let $\Omega$ be a bounded domain in $\mathbf{C}^n$ and $L^2_{\text{hol}}(\Omega)$ the subspace of all holomorphic functions in $L^2(\Omega)$. It is a standard consequence of the mean-value property of holomorphic functions that the point evaluations $f \mapsto f(z)$, $z \in \Omega$, are bounded linear functionals on $L^2_{\text{hol}}(\Omega)$, and thus $L^2_{\text{hol}}(\Omega)$ possesses a *reproducing kernel* — the Bergman kernel $B(x,y)$ of $\Omega$; namely, $B(\cdot, y) \in L^2_{\text{hol}}(\Omega)$ for each $y$, and

$$f(y) = \langle f, B(\cdot, y)\rangle = \int_\Omega f(x) B(y,x)\, dx \qquad \forall f \in L^2_{\text{hol}}(\Omega), \forall y \in \Omega.$$

If $\Omega$ is sufficiently nice (smoothly bounded and strictly pseudoconvex), then the celebrated theorem of Fefferman [20] (with later improvements by Boutet de Monvel and Sjöstrand [10]) gives a description of the boundary singularity of the Bergman kernel: namely, there exist functions $a, b \in C^\infty(\overline{\Omega} \times \overline{\Omega})$ such that

(1) $$B(x,y) = \frac{a(x,y)}{\rho(x,y)^{n+1}} + b(x,y) \log \rho(x,y) \qquad \forall x, y \in \Omega.$$

Here $\rho(x,y) \in C^\infty(\overline{\Omega} \times \overline{\Omega})$ is such that $\partial \rho(x,y)/\partial y$ and $\partial \rho(x,y)/\partial \overline{x}$ vanish to infinite order on the diagonal $x = y$, while $\rho(x,x) = \rho(x)$ is a *defining function* for $\Omega$ in the sense that $\rho(x) > 0$ for $x \in \Omega$ and $\rho(x) = 0$, $\|\nabla \rho(x)\| \neq 0$ for $x \in \partial \Omega$. (It is a consequence of the strict pseudoconvexity of $\Omega$ that such $\rho(x,y)$ exists and can be chosen such that $\text{Re}\, \rho(x,y) > 0$ $\forall x, y \in \Omega$, so that $\log \rho$ can be defined unambiguously.)

1991 *Mathematics Subject Classification.* Primary 32A36; Secondary 32W25, 31B05, 46E22.
Research supported by GA AV ČR grant no. IAA100190802, GA ČR grant no. 201/12/0426 and Czech Ministry of Education research plan no. MSM4781305904.

Fefferman's expansion (1) has subsequently found far-reaching applications in function theory of several complex variables, complex geometry, mathematical physics, operator theory on function spaces, and many other areas (see [4], [21], [25], [27], [12] and [38] for a sample).

The aim of this paper is to give an analogous description for the *harmonic* Bergman kernel, that is, for the reproducing kernel $H(x,y)$ of the subspace

$$L^2_{\text{harm}}(\Omega) := \{f \in L^2(\Omega) : f \text{ is harmonic on } \Omega\}$$

on a bounded domain $\Omega \subset \mathbf{R}^n$, $n \geq 2$, with smooth boundary. Thus $H(\cdot, y) \in L^2_{\text{harm}}(\Omega)$ $\forall y \in \Omega$ and

$$f(y) = \langle f, H(\cdot, y) \rangle = \int_\Omega f(x) H(y,x) \, dx \qquad \forall f \in L^2_{\text{harm}}(\Omega), \forall y \in \Omega.$$

(The existence of the kernel $H(x,y)$ is again a standard, and easily verified, consequence of the mean-value property of harmonic functions.)

Although harmonic Bergman kernels have been around for quite a while — starting probably with the book [5] — and there exist explicit formulas for them in a few special cases (such as the ball and the half-space [15] [2]), as well as for their weighted analogues [34] [30], and various applications of these in operator theory and function theory [14] [35] [13], the description of the boundary behaviour of $H(x,y)$ in the general case seems to be lacking. To the author's knowledge, the only result in this direction in the literature is due to Kang and Koo [31], who gave estimates for the growth of $H(x,y)$ and its derivatives at the boundary: namely, for any multiindices $\alpha, \beta \in \mathbf{N}^n$,

$$\left|\frac{\partial^{|\alpha|+|\beta|} H(x,y)}{\partial x^\alpha \partial y^\beta}\right| \leq \frac{c_{\alpha\beta}}{[d(x) + d(y) + |x-y|]^{n+|\alpha|+|\beta|}}.$$

Here $d(x) = \text{dist}(x, \partial\Omega)$ is the distance of $x$ from the boundary. On the diagonal $x = y$, there are also analogous estimates from below. The main ingredient in their proof is the scaling method, familiar in the theory of several complex variables.

The harmonic Bergman kernel is closely related to another familiar object in analysis, namely, to the Poisson kernel $K(x,\zeta) \equiv K_x(\zeta)$, which gives the solution to the Dirichlet problem

(2) $$\Delta f = 0, \qquad f|_{\partial\Omega} = u$$

on $\Omega$:

(3) $$f(x) = \langle u, K_x \rangle_{\partial\Omega} = \int_{\partial\Omega} u(\zeta) K(x,\zeta) \, d\zeta,$$

where $d\zeta$ denotes the surface measure on $\partial\Omega$. The boundary behaviour of $K(x,\zeta)$ — which, again, is well known in a handful of special cases (like those mentioned before for $H(x,y)$) due to explicit formulas available — has recently been studied by Krantz [32], who showed (using again the scaling method) that

(4) $$K(x,\zeta) \asymp \frac{c_n d(x)}{|x-\zeta|^n}$$

as $x$ and $\zeta$ approach the same point on the boundary; here

$$c_n = \frac{\Gamma(\frac{n}{2})}{\pi^{n/2}}$$

is a constant depending only on the dimension $n$. Again, analogous estimates for the derivatives are also given, as well as numerous references to further information in the literature.

The relationship between $H(x,y)$ and $K(x,\zeta)$ is as follows. Consider the Poisson extension operator from (2)

(5) $$\mathbf{K} : u \mapsto f, \quad \Delta f = 0, \ f|_{\partial\Omega} = u,$$

as an operator from $L^2(\partial\Omega)$ into $L^2(\Omega)$. (In other words, $\mathbf{K}$ is the integral operator with kernel $K$, and $K$ is the Schwartz (distributional) kernel of $\mathbf{K}$.) It can be shown that $\mathbf{K}$ is bounded (even compact), and we denote by $\mathbf{K}^*$ its adjoint. The operator

$$\Lambda := \mathbf{K}^*\mathbf{K}$$

is then an injective compact operator on $L^2(\partial\Omega)$ with dense range, $K_x$ belongs to this range for each $x \in \Omega$, and

(6) $$H(x,y) = \langle \Lambda^{-1} K_x, K_y \rangle_{\partial\Omega} = \mathbf{K}\Lambda^{-1}K_x(y) \qquad \forall x,y \in \Omega.$$

Phrased yet another way, $H(x,y)$ is the distributional (Schwartz) kernel of the operator

(7) $$\mathbf{K}\Lambda^{-1}\mathbf{K}^*,$$

which is nothing but the orthogonal projection (harmonic Bergman projection) $\mathbf{\Pi}_{\mathrm{harm}}$ in $L^2(\Omega)$ onto $L^2_{\mathrm{harm}}(\Omega)$.

Our main results are the following. Let $d$ be a smooth function on $\overline{\Omega}$ such that $d > 0$ on $\Omega$ and $d(x) = \mathrm{dist}(x, \partial\Omega)$ for $x$ near the boundary, and let further $\tilde{x}$ denote the reflection of $x \in \Omega$ with respect to $\partial\Omega$; the latter is well defined for $x$ near the boundary (see §4.3 below for the details). Finally, let $\mathbf{S}^{n-1}$ denote the unit sphere in $\mathbf{R}^n$, $\mathbf{R}_+ = (0, +\infty)$ and $\overline{\mathbf{R}_+} = [0, +\infty)$.

**Theorem 1.** *Let $\Omega \subset \mathbf{R}^n$ be bounded with smooth boundary. Then there exist functions $F \in C^\infty(\partial\Omega \times \overline{\mathbf{R}_+} \times \mathbf{S}^{n-1})$, $G \in C^\infty(\overline{\Omega} \times \partial\Omega)$, with $F(\zeta, 0, \nu) = 1$ for all $\zeta \in \partial\Omega$ and $\nu \in \mathbf{S}^{n-1}$, such that*

(8) $$K(x,\zeta) = \frac{c_n d(x)}{|x-\zeta|^n} \left[ F\Big(\zeta, |x-\zeta|, \frac{x-\zeta}{|x-\zeta|}\Big) + |x-\zeta|^n G(x,\zeta) \log|x-\zeta| \right]$$

*for all $x \in \Omega$, $\zeta \in \partial\Omega$.*

**Theorem 2.** *Let $\Omega \subset \mathbf{R}^n$ be bounded with smooth boundary. Then there exist functions $F \in C^\infty(\overline{\Omega} \times \overline{\Omega} \times \overline{\mathbf{R}_+} \times \mathbf{S}^{n-1})$, $G \in C^\infty(\overline{\Omega} \times \overline{\Omega})$, with $F(x,x,0,\nu) = n\langle \nu, \nabla d(x) \rangle^2 - 1$ for all $x \in \partial\Omega$ and $\nu \in \mathbf{S}^{n-1}$, such that*

(9) $$H(x,y) = \frac{2c_n}{|x-\tilde{y}|^n} F\Big(x, y, |x-\tilde{y}|, \frac{x-\tilde{y}}{|x-\tilde{y}|}\Big) + G(x,y) \log|x-\tilde{y}|$$

*for all $x, y \in \Omega$ close to the boundary.*

Examples are also given showing that, as with the original Fefferman expansion (1) in the holomorphic case, the logarithmic term in (8) and (9) is in general nonzero, even though it is absent in some special situations (e.g. for $n=2$ in Theorem 1).

We remark that the statement of Theorem 1 is equivalent to the same but with $F(x, |x-\zeta|, \frac{x-\zeta}{|x-\zeta|})$, $F \in C^\infty(\overline{\Omega} \times \overline{\mathbf{R}_+} \times \mathbf{S}^{n-1})$, depending on $x$ instead of $\zeta$, or even with $F(x, \zeta, |x-\zeta|, \frac{x-\zeta}{|x-\zeta|})$, $F \in C^\infty(\overline{\Omega} \times \partial\Omega \times \overline{\mathbf{R}_+} \times \mathbf{S}^{n-1})$, allowed to depend on both variables (with $F(\zeta, \zeta, 0, \nu) = 1$ for all $\zeta \in \partial\Omega$ and $\nu \in \mathbf{S}^{n-1}$). Similarly, one can replace the $F(x, y, |x-\tilde{y}|, \frac{x-\tilde{y}}{|x-\tilde{y}|})$ in Theorem 2 by either $F(x, |x-\tilde{y}|, \frac{x-\tilde{y}}{|x-\tilde{y}|})$ or $F(y, |x-\tilde{y}|, \frac{x-\tilde{y}}{|x-\tilde{y}|})$.

Note that the leading order terms in (8) and (9) recover, of course, the coarser estimates from [32] and [31], respectively, recalled above.

As already mentioned, the main tool for deriving the estimates in [32] and [31] was the scaling method; in [32] there is also a rough sketch of a possible way how to obtain (4) — and, in the present author's opinion, also (8) and (9) — by means of a microlocal reduction to a model case using Fourier integral operators (see e.g. §8 in [9], or [10], for a sample of this kind of techniques). We employ a different approach here (avoiding, in particular, the sophisticated machinery of Fourier integral operators), relying on the extension, originating with the work of Boutet de Monvel and others, of the standard calculus of pseudodifferential operators to boundary value problems ("boundary $\Psi$DOs"). This extension includes the "potential" operators like our $\mathbf{K}$, as well as "trace" operators like $\mathbf{K}^*$, $\Psi$DOs on $\partial\Omega$ (which turns out to be the case of the operator $\Lambda = \mathbf{K}^*\mathbf{K}$ above), and, finally, the so-called "singular Green operators", of which an example is precisely the harmonic Bergman projection $\mathbf{\Pi}_{\mathrm{harm}} = \mathbf{K}\Lambda^{-1}\mathbf{K}^*$. The question of identifying the singularities of the distributional kernels of all these operators (including thus, among others, the distributional kernels $K(x, \zeta)$ of $\mathbf{K}$ and $H(x, y)$ of $\mathbf{K}\Lambda^{-1}\mathbf{K}^*$) can then be handled — as in the case of ordinary $\Psi$DOs, where such results go back to Calderón [11] and others — with the aid of the well-known standard apparatus of homogeneous distributions. In particular, in addition to (8) and (9), we in fact obtain a description of the singularities of the distributional kernels of general potential, trace and singular Green operators from Boutet de Monvel's calculus. This includes, by the way, also weighted harmonic Bergman kernels $H_w(x, y)$ on $\Omega$ with respect to weights $w \in C^\infty(\Omega)$ that have "the same order of vanishing" at all points of the boundary (i.e. $w(x) = d(x)^m g(x)$, where $d$ is as above, $m > -1$ and $g \in C^\infty(\overline{\Omega})$ is positive on $\overline{\Omega}$), as well as reproducing kernels of some harmonic Sobolev spaces; see Section 7 below for the details.

It should be noted that all the methods used here are well-established and, in fact, the real-analytic case is more or less treated in [7] (with hints about the $C^\infty$ situation in §7 there); however, the final results seem not to be available anywhere in the literature, nor to be a part of common knowledge, which was the author's reason for writing this paper. Since the paper hopes to aim also at some audience from operator theory and complex function theory, we also decided to include a fairly more extensive review of the various prerequisites than would be strictly necessary for experts in the area of $\Psi$DOs.

The paper is organized as follows. Section 2 reviews the necessary material on boundary $\Psi$DOs (Boutet de Monvel's calculus). Various facts about the Poisson kernel and the harmonic Bergman kernel are collected in Section 3. Singularity resolution of the Schwartz kernel of a general boundary $\Psi$DO is given in Section 4.

Some more specific calculations, which allow us to refine these results to obtain Theorem 1 and Theorem 2, are presented in Section 5 and Section 6, respectively. The final Section 7 contains miscellaneous concluding remarks, open problems, and the like.

NOTATION. Throughout the paper, $\mathbf{Z}$ and $\mathbf{N}$ stand for the integers and the nonnegative integers, respectively. For a multiindex $\alpha \in \mathbf{N}^n$, $\partial_x^\alpha$ (or just $\partial^\alpha$ if the variable is clear from the context) stands for the partial derivative $\partial^{|\alpha|}/\partial x^\alpha = \partial^{|\alpha|}/\partial x_1^{\alpha_1} \ldots \partial x_n^{\alpha_n}$ on $\mathbf{R}^n$, and $D^\alpha = (-i)^{|\alpha|}\partial^\alpha$. The Euclidean norm of $x \in \mathbf{R}^n$ is denoted by $\|x\|$ or simply even $|x|$; the Euclidean scalar product of $x, y \in \mathbf{R}^n$ is denoted by $\langle x, y \rangle$ or $x \cdot y$. The inner products in $L^2(\Omega)$ and $L^2(\partial\Omega)$ are denoted by $\langle \cdot, \cdot \rangle_\Omega$ and $\langle \cdot, \cdot \rangle_{\partial\Omega}$, respectively, with the subscript omitted if there is no danger of confusion. Finally, unless explicitly stated otherwise, $\Omega \subset \mathbf{R}^n$ will be a bounded domain with smooth ($=C^\infty$) boundary.

## 2. Boutet de Monvel calculus

**2.1 Pseudodifferential operators.** Recall that a pseudodifferential operator ($\Psi$DO for short) on $\mathbf{R}^n$ is an operator of the form

$$\tag{10} Au(x) = \int_{\mathbf{R}^n} e^{ix\cdot\xi} a(x,\xi)\,\hat{u}(\xi)\,\dj\xi.$$

Here $\dj\xi = (2\pi)^{-n}d\xi$ is a renormalization of the Lebesgue measure $d\xi$, $x \cdot \xi = x_1\xi_1 + \cdots + x_n\xi_n$, and

$$\tag{11} \hat{u}(\xi) = \int_{\mathbf{R}^n} u(x)\,e^{-ix\cdot\xi}\,dx$$

is the Fourier transform of $u$. The operator (10) is also written as $a(x,D)$, a notation justified by the fact that $A = f(x)D^\alpha$ for $a(x,\xi) = f(x)\xi^\alpha$. The function $a(x,\xi)$, called the *total symbol* of $A$, is usually assumed to lie in Hörmander's class $S^m(\mathbf{R}^n, \mathbf{R}^n)$ for some $m \in \mathbf{R}$, i.e. to belong to $C^\infty(\mathbf{R}^n \times \mathbf{R}^n)$ and satisfy

$$\tag{12} \sup_{x \in K, \xi \in \mathbf{R}^n} \frac{|\partial_x^\alpha \partial_\xi^\beta a(x,\xi)|}{(1+|\xi|)^{m-|\beta|}} < \infty$$

for any compact subset $K$ of $\mathbf{R}^n$ and any multiindices $\alpha, \beta$. It is standard that the operator $a(x,D)$ is then well-defined for any $u$ in, say, the Schwartz space $\mathcal{S}(\mathbf{R}^n)$, and extends by continuity to a bounded operator from the Sobolev space $W^s_{\text{comp}}(\mathbf{R}^n)$ into $W^{s-m}_{\text{loc}}(\mathbf{R}^n)$, for any $s \in \mathbf{R}$, as well as from $\mathcal{D}(\mathbf{R}^n)$ into $C^\infty(\mathbf{R}^n)$ and from the space $\mathcal{D}'_{\text{comp}}(\mathbf{R}^n)$ of compactly supported distributions into $\mathcal{D}'(\mathbf{R}^n)$. Most $\Psi$DOs in this paper will be *classical* (or *polyhomogeneous*), meaning that their total symbol has an asymptotic expansion

$$\tag{13} a(x,\xi) \sim \sum_{j=0}^\infty a_{m-j}(x,\xi),$$

where $a_{m-j}$ is $C^\infty$ in $x, \xi$ and positive homogeneous of degree $m - j$ in $\xi$ for $|\xi| \geq 1$. The notation "$\sim$" means, by definition, that the difference between $a$ and $\sum_{j=0}^{k-1} a_{m-j}$ should belong to the Hörmander class $S^{m-k}$, for each $k = 0, 1, 2, \ldots$. The set of all classical $\Psi$DOs as above (i.e. of order $m$) will be denoted by $\Psi_{\text{cl}}^m$; the (larger) class of all (not necessarily classical) $\Psi$DOs with total symbol in $S^m$ will be denoted by $\Psi^m$, and we set, as usual, $\Psi_{\text{cl}} := \bigcup_{m \in \mathbf{R}} \Psi_{\text{cl}}^m$, $\Psi := \bigcup_{m \in \mathbf{R}} \Psi^m$, and $\Psi^{-\infty} := \bigcap_{m \in \mathbf{R}} \Psi_{\text{cl}}^m = \bigcap_{m \in \mathbf{R}} \Psi^m$. (The classes $\Psi_{\text{cl}}^m$ make sense even for complex $m$, but we will not need those in this paper.) The operators in $\Psi^{-\infty}$ are *smoothing* operators, i.e. their Schwartz kernel is in $C^\infty(\mathbf{R}^n \times \mathbf{R}^n)$.

Combining (10) and (11), we can write the definition of $a(x, D) = A$ as

$$(14) \quad Au(x) = \iint e^{i(x-y)\cdot\xi} a(x, \xi) u(y) \, dy \, d\xi,$$

with the double integral interpreted in suitable sense. One can then define a more general variant of $\Psi$DOs $A = a(x, y, D)$ allowing the symbol $a$ (called then an "*amplitude*") to depend also on $y$:

$$Au(x) = \iint e^{i(x-y)\cdot\xi} a(x, y, \xi) u(y) \, dy \, d\xi,$$

where $a \in C^\infty(\mathbf{R}^{3n})$ is again assumed to satisfy (12) for some $m \in \mathbf{R}$ (with $\partial_x^\alpha$ replaced by $\partial_x^\alpha \partial_y^\gamma$, and the supremum extended also over all $y \in K$). It turns out, nonetheless, that one obtains essentially the same class of operators in this way: namely, for any amplitude $a(x, y, \xi) \in S^m(\mathbf{R}^{2n}, \mathbf{R}^n)$ there exists a symbol $b(x, \xi)$ in $S^m(\mathbf{R}^n, \mathbf{R}^n)$ such that $a(x, y, D) - b(x, D)$ is a smoothing operator. Modulo smoothing operators, $b(x, y)$ is uniquely determined by the asymptotic expansion

$$b(x, \xi) \sim \sum_\alpha \frac{1}{\alpha!} \partial_y^\alpha D_\xi^\alpha a(x, y, \xi)|_{y=x},$$

where the summation extends over all multiindices $\alpha$.

A $\Psi$DO $A = a(x, y, D)$ is *properly supported* if both $A$ and $A^* = \overline{a(y, x, D)}$ (the bar stands for complex conjugation) have the following property: for each compact set $K \subset \mathbf{R}^n$, there exists a compact $K' \subset \mathbf{R}^n$ such that distributions supported in $K$ are mapped into distributions supported in $K'$. Every $\Psi$DO can be written as the sum of a properly supported $\Psi$DO and a smoothing operator.

In general, $\Psi$DOs do not preserve the support of a function (or distribution), however, they have the *pseudolocality* property: if $U \subset \mathbf{R}^n$ is an open set and $u$ is $C^\infty$ on $U$, then $Au$ is also $C^\infty$ on $U$.

The composition of two (classical) $\Psi$DOs, at least one of which is properly supported, is again a (classical) $\Psi$DO: if $a \in \Psi_{\text{cl}}^m$, $b \in \Psi_{\text{cl}}^k$, we have $a(x, D)b(x, D) = c(x, D)$ modulo smoothing operators, where $c \in \Psi_{\text{cl}}^{m+k}$ satisfies

$$(15) \quad c(x, \xi) \sim \sum_\alpha \frac{1}{\alpha!} D_\xi^\alpha a(x, \xi) \partial_x^\alpha b(x, \xi).$$

(The right-hand side determines $c$ only modulo symbols in $S^{-\infty}$.)

$\Psi$DOs behave nicely under coordinate transformations and, consequently, one can define the classes $\Psi^m$, $\Psi_{\text{cl}}^m$ also on manifolds. (Namely, for any partitions of unity $\{\phi_j\}$, $\{\psi_j\}$ subordinate to some atlas $\{\Phi_j\}$ of local coordinate charts, with $\psi_j \equiv 1$ on $\operatorname{supp}\phi_j$, one declares an operator $A$ on the manifold to belong to $\Psi^m$ etc. if and only if the operators

$$\Phi_j^{-1}\psi_j A\phi_k\Phi_k \tag{16}$$

are in $\Psi^m(\mathbf{R}^n)$ etc. for all $j,k$. Note that the operators (16) are automatically properly supported.) In particular, we have the classes $\Psi^m_{\mathrm{cl}}(\partial\Omega)$, $m\in\mathbf{R}$, of classical $\Psi$DOs on the boundary $\partial\Omega$ of a smoothly bounded domain $\Omega$. Note that an operator $A\in\Psi^m_{\mathrm{cl}}(\partial\Omega)$ is then continuous from the Sobolev space $W^s_{\mathrm{comp}}(\partial\Omega)$ into $W^{m-s}_{\mathrm{loc}}(\partial\Omega)$, for any $s\in\mathbf{R}$; in particular, if $\Omega$ is bounded, $A\in\Psi^m_{\mathrm{cl}}(\partial\Omega)$ is continuous from $W^s(\partial\Omega)$ into $W^{m-s}(\partial\Omega)$.

Throughout the rest of this paper, we will almost exclusively be concerned with classical $\Psi$DOs on $\partial\Omega$, and we will thus usually abbreviate $W^s(\partial\Omega)$, $\Psi^m(\partial\Omega)$, etc., just to $W^s$, $\Psi^m$, and so on, and omit the adjective "classical".

Standard references for the material above are e.g. the books by Hörmander [29], Tréves [41], Shubin [40], Folland [23] or Grubb [26]. (All these give also a treatment of $\Psi$DOs on the domain $\Omega\subset\mathbf{R}^n$, rather than on its boundary $\partial\Omega$; this however will not be needed here.)

**2.2 Boundary $\Psi$DOs.** There is an extension of the theory from §2.1 to manifolds with boundary, which we now review; a good reference for the material below and in §§2.3 and 2.4 is Rempel and Schulze [37] and Grubb [26].

Consider first the case of the upper half-space $\mathbf{R}^n_+ = \mathbf{R}^{n-1}\times\mathbf{R}_+$ in $\mathbf{R}^n$, $\mathbf{R}^n_+ = \{x=(x',x_n)\in\mathbf{R}^{n-1}\times\mathbf{R}: x_n>0\}$. Let us denote by $\mathcal{S}_+$ the space of restrictions to $\overline{\mathbf{R}_+}=[0,+\infty)$ of functions in the Schwartz space $\mathcal{S}(\mathbf{R})$, and, similarly, by $\mathcal{S}_{++}$ the space of restrictions to $\overline{\mathbf{R}_+}\times\overline{\mathbf{R}_+}$ of functions in the Schwartz space $\mathcal{S}(\mathbf{R}^2)$. For $d\in\mathbf{R}$, the space $S^d(\mathbf{R}^m,\mathbf{R}^{n-1},\mathcal{S}_+)$ consists, by definition, of all $C^\infty$ functions $f$ of $(X,x_n,\xi')\in\mathbf{R}^m\times\overline{\mathbf{R}_+}\times\mathbf{R}^{n-1}$, lying in $\mathcal{S}_+$ with respect to $x_n$, such that for all $\alpha\in\mathbf{N}^{n-1}$, $\beta\in\mathbf{N}^m$, $k,l\in\mathbf{N}$ and compact subsets $K$ of $\mathbf{R}^m$,

$$\sup_{X\in K,\, x_n\geq 0,\, \xi'\in\mathbf{R}^{n-1}} \frac{|x_n^k\partial_{x_n}^l\partial_X^\beta\partial_{\xi'}^\alpha f(X,x_n,\xi')|}{(1+|\xi'|)^{d+1+l-k-|\alpha|}} < \infty.$$

The space $S^d(\mathbf{R}^m,\mathbf{R}^{n-1},\mathcal{S}_{++})$ consists similarly, by definition, of all $C^\infty$ functions $g$ of $(X,x_n,y_n,\xi')\in\mathbf{R}^m\times\overline{\mathbf{R}_+}\times\overline{\mathbf{R}_+}\times\mathbf{R}^{n-1}$, lying in $\mathcal{S}_{++}$ with respect to $(x_n,y_n)$, such that for all $\alpha\in\mathbf{N}^{n-1}$, $\beta\in\mathbf{N}^m$, $k,l,p,q\in\mathbf{N}$ and compact subsets $K$ of $\mathbf{R}^m$,

$$\sup_{X\in K,\, x_n,y_n\geq 0,\, \xi'\in\mathbf{R}^{n-1}} \frac{|x_n^k\partial_{x_n}^l y_n^p\partial_{y_n}^q\partial_X^\beta\partial_{\xi'}^\alpha g(X,x_n,y_n,\xi')|}{(1+|\xi'|)^{d+2+l+q-k-p-|\alpha|}} < \infty.$$

The subspaces $S^d_{\mathrm{cl}}$ of classical (or polyhomogeneous) elements in these $S^d$ are defined to consist, analogously as for symbols of $\Psi$DOs in §2.1, of all $f,g$ with asymptotic expansions

$$f\sim\sum_{l=0}^\infty f_{d-l} \quad\text{or}\quad g\sim\sum_{l=0}^\infty g_{d-l},$$

respectively, where

$$f_j(X,\tfrac{x_n}{\lambda},\lambda\xi') = \lambda^{j+1}f(X,x_n,\xi'), \qquad g_j(X,\tfrac{x_n}{\lambda},\tfrac{y_n}{\lambda},\lambda\xi') = \lambda^{j+2}g(X,x_n,y_n,\xi'),$$

for $\lambda \geq 1$ and $|\xi'| \geq 1$.

A *negligible potential operator* is an operator from functions on $\mathbf{R}^{n-1}$ into functions on $\mathbf{R}^n_+$ whose distributional (Schwartz) kernel is smooth on $\overline{\mathbf{R}^n_+} \times \mathbf{R}^{n-1}$. A *potential operator* of order $d \in \mathbf{R}$ is an operator from functions on $\mathbf{R}^{n-1}$ into functions on $\mathbf{R}^n_+$ which is a sum of a negligible potential operator and an operator defined by

$$(17) \qquad Ku(x', x_n) = \int_{\mathbf{R}^{n-1}} e^{ix'\cdot\xi'} k(x', x_n, \xi') \, \hat{u}(\xi') \, d\xi'$$

where $k \in S^{d-1}_{\mathrm{cl}}(\mathbf{R}^{n-1}, \mathbf{R}^{n-1}, \mathcal{S}_+)$. As with $\Psi$DO, one gets the same class of operators using kernels $k(y', x_n, \xi')$ instead of $k(x', x_n, \xi')$, and even symbols (amplitudes) depending on both $x'$ and $y'$ can be allowed:

$$Ku(x', x_n) = \int_{\mathbf{R}^{n-1}} \int_{\mathbf{R}^{n-1}} e^{i(x'-y')\cdot\xi'} k(x', y', x_n, \xi') \, u(y') \, dy' \, d\xi',$$

with $k \in S^{d-1}_{\mathrm{cl}}(\mathbf{R}^{2n-2}, \mathbf{R}^{n-1}, \mathcal{S}_+)$ and the integral suitably interpreted (as an oscillatory integral), but again yield the same class of operators.

A *negligible trace operator* of class $r \in \mathbf{N}$ is an operator from functions on $\mathbf{R}^n_+$ into functions on $\mathbf{R}^{n-1}$ of the form

$$\sum_{0 \leq j < r} s_j \gamma_j + T,$$

where $\gamma_j$ is the $j$-th order boundary normal derivative $\gamma_j u(x') = (\partial^j_{x_n} u)(x', 0)$, and $s_j, T$ have distributional kernels that are smooth on $\mathbf{R}^{n-1} \times \mathbf{R}^{n-1}$ and $\mathbf{R}^{n-1} \times \overline{\mathbf{R}^n_+}$, respectively. A *trace operator* of class $r \in \mathbf{N}$ and order $d \in \mathbf{R}$ is an operator from functions on $\mathbf{R}^n_+$ into functions on $\mathbf{R}^{n-1}$ which is a sum of a negligible trace operator and an operator of the form

$$(18) \qquad \sum_{0 \leq j < r} s_j \gamma_j + T,$$

where $s_j \in \Psi^{d-j}_{\mathrm{cl}}(\mathbf{R}^{n-1})$ and

$$Tu(x') = \int_{\mathbf{R}^{n-1}} \int_0^\infty e^{ix'\cdot\xi'} t(x', x_n, \xi') \, \acute{u}(\xi', x_n) \, dx_n d\xi',$$

where $\acute{u}(\xi', x_n)$ stands for the partial Fourier transform of $u(x', x_n)$ with respect to the $x'$ variable, and $t \in S^d_{\mathrm{cl}}(\mathbf{R}^{n-1}, \mathbf{R}^{n-1}, \mathcal{S}_+)$. Again, symbols $t(y', x_n, \xi')$ or more general amplitudes $t(x', y', x_n, \xi')$ from $S^d_{\mathrm{cl}}(\mathbf{R}^{2n-2}, \mathbf{R}^{n-1}, \mathcal{S}_+)$ can also be allowed, but lead to the same class of operators.

A *negligible singular Green operator* of class $r \in \mathbf{N}$ is an operator on $\mathbf{R}^n_+$ of the form

$$\sum_{0 \leq j < r} K_j \gamma_j + G,$$

where $K_j$ are negligible potential operators and the distributional kernel of $G$ is smooth on $\overline{\mathbf{R}^n_+} \times \overline{\mathbf{R}^n_+}$. A *singular Green operator* of class $r \in \mathbf{N}$ and order $d \in \mathbf{R}$ is an operator on $\mathbf{R}^n_+$ which is a sum of a negligible singular Green operator and an operator of the form

$$\sum_{0 \leq j < r} K_j \gamma_j + G, \tag{19}$$

where $K_j$ are potential operators of order $d - j$ and

$$Gu(x', x_n) = \int_{\mathbf{R}^{n-1}} \int_0^\infty e^{ix' \cdot \xi'} g(x', x_n, y_n, \xi') \, \acute{u}(\xi', y_n) \, dy_n \, đ\xi', \tag{20}$$

where $g \in S_{\mathrm{cl}}^{d-1}(\mathbf{R}^{n-1}, \mathbf{R}^{n-1}, \mathcal{S}_{++})$. Again, symbols $g(y', x_n, y_n, \xi')$ or more general amplitudes $g(x', y', x_n, y_n, \xi')$ from $S_{\mathrm{cl}}^{d-1}(\mathbf{R}^{n-1} \times \mathbf{R}^{n-1}, \mathbf{R}^{n-1}, \mathcal{S}_{++})$ can be allowed but lead to the same class of operators.

We denote the three types of operators just introduced by $\mathcal{K}^d(\overline{\mathbf{R}^n_+})$, $\mathcal{T}^d_r(\overline{\mathbf{R}^n_+})$, and $\mathcal{G}^d_r(\overline{\mathbf{R}^n_+})$, respectively.

As with $\Psi$DOs, the three types of operators above are pseudolocal in the $x'$-variable, in the obvious sense. They are called "*boundary $\Psi$DOs*" (b$\Psi$DOs).

We will call the functions $k, t, g$ the *symbols* of $K$, $T$ and $G$, respectively.[1]

The calculus of b$\Psi$DOs, initiated by Boutet de Monvel [6], [8], implies, first of all, that the three classes of operators defined above behave nicely under compositions and taking adjoints; and, second, that solution operators to boundary value problems — like our Poisson operator $\mathbf{K}$ which we had before — belong to this calculus. Let us briefly describe only those details about this that we will need.

*Remark.* There exist also more general "nonclassical" versions of $T$, $K$ and $G$, but they will not be needed in this paper. □

**2.3 Composition rules.** Let $K$ be a potential operator, $T$ a trace operator of class zero, and $S$ a (classical) $\Psi$DO on $\mathbf{R}^{n-1}$. Then

$$\begin{aligned} & KS \text{ is a potential operator,} \\ & ST \text{ is a trace operator,} \\ & TK \text{ is a } \Psi\text{DO on } \partial\Omega, \text{ and} \\ & KT \text{ is a singular Green operator,} \end{aligned} \tag{21}$$

with orders adding up, provided at least one of the factors in each case is properly supported in the $x'$-variable. The symbols of these products are given by the usual formula (15) in the $x'$-variable, but new things happen in the $x_n$-variable. Let us define

$$\begin{aligned} (k \circ_n t)(x', x_n, y_n, \xi') &:= k(x', x_n, \xi') t(x', y_n, \xi'), \\ (t \circ_n k)(x', \xi') &:= \int_0^\infty t(x', x_n, \xi') k(x', x_n, \xi') \, dx_n, \\ (k \circ_n s)(x', x_n, \xi') &:= k(x', x_n, \xi') s(x', \xi'), \\ (s \circ_n t)(x', x_n, \xi') &:= s(x', \xi') t(x', x_n, \xi'). \end{aligned} \tag{22}$$

---

[1] The reader is warned that this is somewhat at odds with the standard notations and terminology in the area, including the references such as [26] mentioned further down: namely, our $k$ is usually denoted by $\tilde{k}$ and called (for good reasons) *symbol-kernel* rather than symbol, of $K \in \mathcal{K}^d$; the term "symbol" being reserved for the Fourier transform of our $k(x', x_n, \xi')$ with respect to the $x_n$ variable, after extending by zero for $x_n < 0$. Similarly for $\mathcal{T}$ and $\mathcal{G}$. Also "potential operators" are called "Poisson operators" in [26] and elsewhere; we reserve here the latter term for the Poisson operator $\mathbf{K}$ from (5).

Let $(A,B)$ be one of the pairs of operators $(K,S)$, $(S,T)$, $(T,K)$ or $(K,T)$, and let $a,b$ be the corresponding symbols $k,t,s$ or $g$. Then $AB = C$ is an operator of the type indicated in (21), with symbol

$$(23) \qquad c \sim \sum_\alpha \frac{1}{\alpha!}(D_{\xi'}^\alpha a) \circ_n (\partial_{x'}^\alpha b).$$

There are also composition formulas for trace operators of class $r > 0$; we will only need the case of $TK$ with $r = 1$ here, which reduces to the simple rule

$$(24) \qquad \gamma_0 K = S \quad \text{with} \quad s(x',\xi') = k(x',0,\xi').$$

Similarly for $\gamma_0 G$ with a singular Green operator $G$. Finally, there are rules for the products $PK = K'$ and $PG = G'$, where $P$ is a differential operator on $\mathbf{R}^n$: namely, on the level of the $x'$ variables one uses the standard $\Psi$DO rule (15), while for $P = x_n^k \partial_{x_n}^l$ one just has

$$k'(x',x_n,\xi') = x_n^k \partial_{x_n}^l k(x',x_n,\xi')$$

and similarly for $g$ and $g'$.

Similar rules hold also for other possible compositions, but again will not be needed here.

In addition to composition, the calculus behaves well with respect to adjoints: namely, the adjoint of a Poisson operator $K \in \mathcal{K}^d$ with symbol $k$ is the trace operator $T \in \mathcal{T}_0^{d-1}$ with symbol

$$(25) \qquad t(x',x_n,\xi') \sim \sum_\alpha \frac{1}{\alpha!} D_{\xi'}^\alpha \partial_{x'}^\alpha \overline{k(x',x_n,\xi')}$$

(the bar denoting complex conjugation). There is also an analogue for $\mathcal{G}_0^d$, which will not be needed here.

**2.4 Boundary value problems.** For $T \in \mathcal{T}_r^m(\overline{\mathbf{R}_+^n})$, the *boundary symbol* of $T$ is the operator $t(x',\xi',D_n)$ from $\mathcal{S}_+$ into $\mathbf{C}$ defined by

$$t(x',\xi',D_n)u = \int_0^\infty t(x',x_n,\xi')\,u(x_n)\,dx_n + \sum_{0 \le j < r} s_j(x',\xi')\gamma_j u.$$

Similarly, the boundary symbol of $G \in \mathcal{G}_r^d(\overline{\mathbf{R}_+^n})$ is the operator on $\mathcal{S}_+$ defined by

$$g(x',\xi',D_n)u(x_n) = \int_0^\infty g(x',x_n,y_n,\xi')\,u(y_n)\,dy_n + \sum_{0 \le j < r} k_j(x',x_n,\xi')\gamma_j u,$$

and the boundary symbol of $K \in \mathcal{K}^q(\overline{\mathbf{R}_+^n})$ is the operator from $\mathbf{C}$ into $\mathcal{S}_+$ defined by

$$k(x',\xi',D_n)a = k(x',x_n,\xi')a, \qquad a \in \mathbf{C}.$$

The *principal boundary symbol* is defined similarly using only the leading-order term (the principal symbol) $t_m$, $g_d$ or $k_q$ in the asymptotic expansion of $t$, $g$ or $k$, respectively.

Now let $A$ be the matrix of operators

$$
(26) \qquad A = \begin{bmatrix} P_+ + G & K \\ T & S \end{bmatrix},
$$

with $K, T, G$ as above, $S$ a $\Psi$DO on $\mathbf{R}^{n-1}$ of order $m+q-d$, and $P_+ = r_+ P e_+$ where $P$ is a differential operator on $\mathbf{R}^n$, of the same order $d$ as $G$, with smooth coefficients, $r_+$ stands for the restriction from $\mathbf{R}^n$ to $\mathbf{R}^n_+$, and $e_+$ for the operator of extension by zero from $\mathbf{R}^n_+$ to $\mathbf{R}^n$; abusing notation slightly, we will sometimes write just $P$ instead of $P_+$. Also, more generally we even allow $K$ to be an $M$-tuple (row matrix) of potential operators, $T$ an $M'$-tuple (column matrix) of trace operators, and $S$ an $M' \times M$ matrix of $\Psi$DOs ($M, M' = 0, 1, 2, \dots$). The *principal interior symbol* of $A$ is just the principal symbol $p_d(x, \xi)$ of $P$; the *principal boundary symbol* is the matrix of operators

$$
\mathfrak{a}(x', \xi', D_n) := \begin{bmatrix} p_d(x', 0, \xi', D_n) + g_d(x', \xi', D_n) & k_q(x', \xi', D_n) \\ t_m(x', \xi', D_n) & s_{m+q-d}(x', \xi'), \end{bmatrix}
$$

with $p_d(x', 0, \xi', D_n) : \mathcal{S}_+ \to \mathcal{S}_+$ defined by

$$
p_d(x', 0, \xi', D_n) u(x_n) = \int_0^\infty e^{i x_n \xi_n} p_d(x', 0, \xi', \xi_n) \widehat{e_+ u}(\xi_n) \, d\xi_n,
$$

where $e_+$ is the operator of "extension by zero" from $\mathbf{R}_+$ to $\mathbf{R}$. The matrix $A$ is said to be *elliptic* (of orders $d, m, q$ and class $r$) if $p_d(x, \xi) \neq 0$ for all $|\xi| \geq 1$ and all $x \in \mathbf{R}^n_+$, and $\mathfrak{a}(x', \xi', D_n) : \mathcal{S}_+ \times \mathbf{C}^M \to \mathcal{S}_+ \times \mathbf{C}^{M'}$ is bijective for all $|\xi'| \geq 1$ and $x' \in \mathbf{R}^{n-1}$.

Finally, the three classes of operators $\mathcal{K}, \mathcal{T}, \mathcal{G}$ are invariant under coordinate transformations of $\overline{\mathbf{R}^n_+}$ preserving the boundary $\partial \overline{\mathbf{R}^n_+} = \mathbf{R}^{n-1} \times \{0\}$. Thanks to this invariance, one can again — via local coordinate charts — define the analogous classes $\mathcal{K}^d$, $\mathcal{T}^d_r$ and $\mathcal{G}^d_r$ for $\overline{\mathbf{R}^n_+}$ replaced by any manifold with boundary; in particular, the classes $\mathcal{K}^d(\overline{\Omega})$, $\mathcal{T}^d_r(\overline{\Omega})$ and $\mathcal{G}^d_r(\overline{\Omega})$ for $\mathbf{R}^n_+$ replaced by a domain $\Omega \subset \mathbf{R}^n$, and the boundary $\partial \overline{\mathbf{R}^n_+} = \mathbf{R}^{n-1}$ replaced by the boundary $\partial \Omega$ of $\Omega$. The definitions of principal boundary symbols and ellipticity carry over to this setting.

With all these definitions, it is then true that

(27) if $A$ is elliptic of orders $d, m, q$ and class $r$ and $\Omega$ is bounded, then there exists a matrix $B$ of operators of the form (26), elliptic of orders $-d, -q, -m$ and class $r' = \max\{r - d, 0\}$, which is a parametrix for $A$, i.e. $AB - I$ and $BA - I$ are negligible.

A potential operator $K \in \mathcal{K}^d(\overline{\Omega})$ whose amplitude is compactly supported in both $x'$ and $y'$ can be shown to map $W^s(\partial \Omega)$ continuously into $W^{s-d+\frac{1}{2}}(\Omega)$, for any $s \in \mathbf{R}$. In particular, if $\Omega$ is bounded, then any $K \in \mathcal{K}^d$ is continuous from $W^s(\partial \Omega)$ into $W^{s-d+1/2}(\Omega)$, $s \in \mathbf{R}$. Similarly, $T \in \mathcal{T}^d_r(\overline{\Omega})$ and $G \in \mathcal{G}^d_r(\overline{\Omega})$ are continuous $W^s(\Omega) \to W^{s-d-1/2}(\partial \Omega)$ and $W^s(\Omega) \to W^{s-d}(\Omega)$, respectively, for any $s > r - \frac{1}{2}$, $s \geq 0$, if $\Omega$ is bounded.

*Remark.* The result (27) actually holds also for suitable pseudodifferential (i.e. not necessarily differential) operators $P$, namely, for those that satisfy the *transmission condition*; but again we will not need this here. $\square$

*Remark.* For $\Omega$ unbounded, (27) gets more complicated: for instance, even for the matrix (26) corresponding to the ordinary Dirichlet problem for $\Delta$ on $\mathbf{R}^n_+$ (cf. the beginning of Section 3), the kernel $k(x', x_n, \xi')$ in (17) of the corresponding Poisson operator $\mathbf{K}$ equals $k(x', x_n, \xi') = e^{-x_n|\xi'|}$, which is not smooth at $\xi' = 0$. Nonetheless, the localized operator $\mathbf{K}\chi$ (i.e. $\mathbf{K}$ preceded by multiplication by $\chi$), with any $\chi \in \mathcal{D}(\mathbf{R}^{n-1})$, will already be an honest potential operator, i.e. with smooth $k$. $\square$

## 3. Some simple facts

In this section, we establish various useful facts about the Poisson operator $\mathbf{K}$, the Poisson kernel $K(x,\zeta)$, the harmonic Bergman projection $\Pi_{\mathrm{harm}}$, and the harmonic Bergman kernel $H(x,y)$, that will be needed later.

Consider first our Poisson (harmonic extension) operator $\mathbf{K}: u \mapsto f$ solving the classical harmonic extension problem

$$\Delta f = 0 \text{ on } \Omega, \quad f = u \text{ on } \partial\Omega;$$

that is,

$$\begin{bmatrix} \Delta \\ \gamma_0 \end{bmatrix} f = \begin{bmatrix} 0 \\ u \end{bmatrix}.$$

The operator

$$A = \begin{bmatrix} \Delta \\ \gamma_0 \end{bmatrix} : C^\infty(\overline{\Omega}) \to \begin{matrix} C^\infty(\overline{\Omega}) \\ \oplus \\ C^\infty(\partial\Omega) \end{matrix}$$

is of the type (26) discussed in §2.4, with interior symbol $p_2(x,\xi) = -|\xi|^2$ and boundary symbol given by

$$\mathfrak{a}(x', \xi', D_n) = \begin{bmatrix} -|\xi'|^2 + \partial^2_{x_n} \\ \gamma_0 \end{bmatrix} : \mathcal{S}_+ \to \begin{matrix} \mathcal{S}_+ \\ \oplus \\ \mathbf{C} \end{matrix}.$$

The equation $\mathfrak{a}(x', \xi', D_n)u = \begin{bmatrix} 0 \\ a \end{bmatrix}$ has the unique solution $u(x_n) = ae^{-x_n|\xi'|}$ in $\mathcal{S}_+$, and similarly $\mathfrak{a}(x', \xi', D_n)u = \begin{bmatrix} v \\ 0 \end{bmatrix}$ is uniquely solvable. Thus $A$ is elliptic and has a parametrix, unique up to negligible operators, of the form $[P + G, K]$, with $P$ a $\Psi$DO on $\mathbf{R}^n$ of order $-2$, $G$ a singular Green operator of order $-2$ and class $\max\{1-2,0\} = 0$, and $K$ a potential operator on $\Omega$ of order $0$. On the other hand, we know that

$$B = [\mathbf{G}, \mathbf{K}],$$

with $\mathbf{K}$ as above and $\mathbf{G}$ the solution operator to the classical Dirichlet problem

$$\Delta \mathbf{G} f = f \text{ on } \Omega, \quad \mathbf{G} f = 0 \text{ on } \partial\Omega,$$

is a right inverse for $A$. It follows that we can take $P + G = \mathbf{G}$ and $\mathbf{K} = K$; that is,

(28) $$\mathbf{K} \in \mathcal{K}^0(\overline{\Omega}).$$

By the facts reviewed in Section 2, it follows that

(29) $$\mathbf{K}: W^s(\partial\Omega) \to W^{s+1/2}(\Omega), \qquad \forall s \in \mathbf{R};$$

in particular, $\mathbf{K}$ is a continuous operator from $W^0(\partial\Omega) = L^2(\partial\Omega)$ into $W^{1/2}_{\text{harm}}(\Omega) \subset L^2_{\text{harm}}(\Omega) \subset L^2(\Omega)$ (even compact). Let $\mathbf{K}^* : L^2(\Omega) \to L^2(\partial\Omega)$ be its Hilbert space adjoint and consider the operator

$$\Lambda := \mathbf{K}^*\mathbf{K}.$$

Note that $\Lambda$ is injective (since $\mathbf{K}$ is). From

(30) $$\Lambda^{-1}\mathbf{K}^*\mathbf{K} = I$$

it is immediate that $\mathbf{K}\Lambda^{-1}\mathbf{K}^*$ is selfadjoint, idempotent, and vanishes on $\operatorname{Ker} \mathbf{K}^* = (\overline{\operatorname{Ran} \mathbf{K}})^\perp$ while being the identity on the range of $\mathbf{K}$. Consequently,

(31) $$\mathbf{K}\Lambda^{-1}\mathbf{K}^* = \mathbf{\Pi}_{\text{harm}}$$

is the orthogonal projection in $L^2(\Omega)$ onto the closure of $\operatorname{Ran} \mathbf{K}$, i.e. the harmonic Bergman projection of $L^2(\Omega)$ onto $L^2_{\text{harm}}(\Omega)$.

*Remark.* Comparing (30) with the definition of $\mathbf{K}$, we also see that

$$\Lambda^{-1}\mathbf{K}^*|_{\operatorname{Ran} \mathbf{K}} =: \gamma$$

coincides with the operator of "taking the nontangential boundary values" of harmonic functions. By elliptic regularity (or Boutet de Monvel's calculus), $\gamma$ again extends to a continuous operator from $W^s_{\text{harm}}(\Omega)$ into $W^{s-1/2}(\partial\Omega)$, for any $s \in \mathbf{R}$, which is the left inverse of $\mathbf{K}$. □

By Boutet de Monvel's calculus, we have from (28) $\mathbf{K}^* \in \mathcal{T}_0^{-1}(\overline{\Omega})$ and thus

$$\Lambda = \mathbf{K}^*\mathbf{K} \in \Psi^{-1}_{\text{cl}}(\partial\Omega).$$

Furthermore, by (22), (23) and (25), the leading symbol $s_{-1}$ of $\Lambda$ is given by, in any local coordinate chart,

$$s_{-1}(x',\xi') = \int_0^\infty |k_0(x', x_n, \xi')|^2 \, dx_n > 0,$$

i.e. $\Lambda$ is elliptic. Its inverse $\Lambda^{-1}$ is thus a positive injective selfadjoint operator in $\Psi^1_{\text{cl}}(\partial\Omega)$. Finally, by Boutet de Monvel's calculus again and (31),

$$\mathbf{\Pi}_{\text{harm}} = \mathbf{K}\Lambda^{-1}\mathbf{K}^* \in \mathcal{G}^0_0(\overline{\Omega})$$

is a singular Green operator on $\Omega$ of class zero and order 0.

We conclude this section by proving one more formula relating the Poisson and the harmonic Bergman kernel.

**Proposition 3.** *For $x, y \in \Omega$,*
$$H(x,y) = \langle \Lambda^{-1} K_x, K_y \rangle_{\partial \Omega}.$$

*Proof.* For any $f \in L^2(\Omega)$ and $x \in \Omega$, we have
$$\mathbf{\Pi}_{\text{harm}} f(x) = \int_\Omega f(y) H(x,y) \, dy.$$

Thus
$$\begin{aligned}
\langle f, H_x \rangle_\Omega &= \mathbf{\Pi}_{\text{harm}} f(x) = \mathbf{K} \Lambda^{-1} \mathbf{K}^* f(x) \\
&= \langle \Lambda^{-1} \mathbf{K}^* f, K_x \rangle_{\partial \Omega} \quad \text{by (3)} \\
&= \langle f, \mathbf{K} \Lambda^{-1} K_x \rangle_\Omega.
\end{aligned}$$

This means that $H_x = \mathbf{K} \Lambda^{-1} K_x$ and
$$H_x(y) = \mathbf{K} \Lambda^{-1} K_x(y) = \langle \Lambda^{-1} K_x, K_y \rangle_{\partial \Omega}$$

as claimed. □

## 4. Schwartz kernels of boundary $\Psi$DOs

**4.1 Homogeneous distributions.** Recall that a function $u \in C^\infty(\mathbf{R}^n \setminus \{0\})$ is said to be *homogeneous* of degree $s$, $s \in \mathbf{R}$, if $u(\lambda x) = \lambda^s u(x)$ for all $x \in \mathbf{R}^n \setminus \{0\}$ and $\lambda > 0$. The definition extends also to distributions on $\mathbf{R}^n \setminus \{0\}$ or $\mathbf{R}^n$ in an obvious way. It is standard (see e.g. [28], Chapter 3, §2) that $u$ can always be prolonged from $\mathbf{R}^n \setminus \{0\}$ to a tempered distribution $\dot{u}$ on $\mathbf{R}^n$ and, furthermore, $\dot{u}$ is uniquely determined and homogeneous if $-s - n \notin \mathbf{N}$, while for $s = -n - k$, $k \in \mathbf{N}$, $\dot{u}$ is not unique (one can add any derivative of order $k$ of the Dirac distribution $\delta$ at the origin) and is not homogeneous in general, but contains a logarithmic term. The Fourier transform of $\dot{u}$ has the form

(32) $$\widehat{\dot{u}} = U - Q \log | \cdot |,$$

where $U$, $U \in C^\infty(\mathbf{R}^n \setminus \{0\})$ if $u \in C^\infty(\mathbf{R}^n \setminus \{0\})$, is a homogeneous distribution on $\mathbf{R}^n$ of degree $-n - s$, while $Q$ is a polynomial of degree $-n - s =: k$ (thus $Q = 0$ if $k \notin \mathbf{N}$) given by

(33) $$Q(\xi) = \int_{\mathbf{S}^{n-1}} \frac{(-ix \cdot \xi)^k u(x)}{k!} \, dx,$$

where $dx$ stands for the unnormalized surface measure on the unit sphere $\mathbf{S}^{n-1}$.

These facts can be neatly summarized as follows: for $s \in \mathbf{R}$, let $\mathcal{H}_s$ be the vector space of all distributions on $\mathbf{R}^n$ of the form
$$\dot{U}_s + Q_s \log | \cdot | + P_{-n-s}(D) \delta_0,$$

where $U_s \in C^\infty(\mathbf{R}^n \setminus \{0\})$ is homogeneous of degree $s$, $Q_s$ is a polynomial of degree $s$ (thus $Q_s = 0$ if $s$ is not a nonnegative integer), $P_{-n-s}$ is a polynomial of degree $-n-s$ (thus $P_{-n-s} = 0$ if $-n-s$ is not a nonnegative integer), and $\delta_0$ stands for the Dirac distribution (unit point mass) at the origin. Then the prolongation $\dot{u}$ of a function $u \in C^\infty(\mathbf{R}^n \setminus \{0\})$ homogeneous of degree $s$ always exists, is unique if $-n-s \notin \mathbf{N}$, and all such prolongations belong to $\mathcal{H}_s$. Noting that the Fourier transform of $\log |x|$ coincides (up to a constant factor) with $|\xi|^{-n}$, it further follows that in fact

(34)    the Fourier transform is a bijection of $\mathcal{H}_s$ onto $\mathcal{H}_{-n-s}$,    $\forall s \in \mathbf{R}$.

The following assertions seem to be standard knowledge, but difficult to pinpoint in the literature in complete generality; see e.g. Neri [36], Theorem 1.5, or Calderon [11], Theorem 28.

Recall that a function $\chi \in \mathcal{D}(\mathbf{R}^n)$, $0 \leq \chi \leq 1$, is called a *cutoff* function if $\chi \equiv 1$ in some neighbourhood of the origin; and $\theta$ is a *patch* function if $\theta = 1 - \chi$, with $\chi$ a cutoff function.

**Proposition 4.** (i) *Let $u \in C^\infty(\mathbf{R}^n \setminus \{0\})$ be homogeneous of degree $s$, $s \in \mathbf{R}$, and let $\theta$ be a patch function. Then*

$$\widehat{\theta u} = (U - Q \log |\cdot|)\chi + g, \tag{35}$$

*where $U$ and $Q$ are as in (32) and (33), $\chi$ is a cutoff function, and $g \in \mathcal{S}(\mathbf{R}^n)$.*

(ii) *Let $\theta, \chi$ be as above. For $s \in \mathbf{R}$, denote*

$$\theta \mathcal{H}_s + \mathcal{S} := \{\theta u + g : u \in \mathcal{H}_s, g \in \mathcal{S}\},$$
$$\chi \mathcal{H}_s + \mathcal{S} := \{\chi u + g : u \in \mathcal{H}_s, g \in \mathcal{S}\}$$

(*these sets are independent of the choice of $\chi$ and $\theta$*). *Then the Fourier transform is a bijection of $\theta \mathcal{H}_s + \mathcal{S}$ onto $\chi \mathcal{H}_{-n-s} + \mathcal{S}$.*

*Proof.* (i) Let $\dot{u}$ be as above; then $\theta \dot{u} = \theta u$ differs from $\dot{u}$ by a distribution with compact support. By the Paley-Wiener theorem, we thus have

$$\widehat{\theta u} = \widehat{\theta \dot{u}} = \hat{\dot{u}} + h,$$

with $h$ an entire function on $\mathbf{C}^n$. Hence, for any cutoff function $\chi$, $\widehat{\theta u} - \chi \hat{\dot{u}} \in C^\infty(\mathbf{R}^n)$.

On the other hand, for any multiindices $\alpha, \beta$, $\xi^\alpha D^\beta \widehat{\theta u} = (-1)^{|\beta|} \widehat{D^\alpha x^\beta \theta u}$ is the Fourier transform of a smooth function on $\mathbf{R}^n$ homogeneous of degree $s + |\beta| - |\alpha|$ near infinity, and, hence, integrable for $|\alpha| > s + |\beta| + n$. Thus by the Riemann-Lebesgue lemma, $\xi^\alpha D^\beta \widehat{\theta u} \in C_0(\mathbf{R}^n)$ for $|\alpha| > s + |\beta| + n$. It follows that $\partial^\beta \widehat{\theta u}$ — and, hence, $\partial^\beta(\widehat{\theta u} - \chi \hat{\dot{u}})$ — decays rapidly at infinity, for any $\beta$. Thus $\widehat{\theta u} - \chi \hat{\dot{u}} \in \mathcal{S}(\mathbf{R}^n)$, proving the first claim.

(ii) Just as in the proof of part (i), one sees that $\widehat{\theta u} - \chi \hat{\dot{u}} \in \mathcal{S}$ for $u \in \mathcal{H}_s$ (the only difference is that $D^\alpha x^\beta \theta u$ need now no longer be homogeneous of degree $s + |\beta| - |\alpha|$ near infinity, but may also contain a term of this form multiplied by $\log |x|$, which however has no effect on its integrability). Thus by (32), $\widehat{\theta \mathcal{H}_s} \subset \chi \mathcal{H}_{-n-s} + \mathcal{S}$, and hence also $[\theta \mathcal{H}_s + \mathcal{S}]^\wedge \subset \chi \mathcal{H}_{-n-s} + \mathcal{S}$.

Conversely, if $v \in \mathcal{H}_{-n-s}$, then by (34) there is $u \in \mathcal{H}_s$ with $v = \hat{u}$, so as we have observed at the beginning of the preceding paragraph, $\widehat{\theta u} - \chi v \in \mathcal{S}$; thus $\chi v \in \widehat{\theta \mathcal{H}_s} + \mathcal{S}$. Hence $\chi \mathcal{H}_{-n-s} + \mathcal{S} \subset \widehat{\theta \mathcal{H}_s} + \mathcal{S} = [\theta \mathcal{H}_s + \mathcal{S}]^\wedge$. □

**Corollary 5.** *If $A \in \Psi^m_{\mathrm{cl}}(\mathbf{R}^n)$ is properly supported, then its Schwartz kernel $k_A$ satisfies for $x \neq y$*

$$k_A(x,y) = |x-y|^{-n-m} F\Big(x, |x-y|, \frac{x-y}{|x-y|}\Big) + G(x,y) \log|x-y|,$$

*where $F \in C^\infty(\mathbf{R}^n \times \overline{\mathbf{R}_+} \times \mathbf{S}^{n-1})$, $G \in C^\infty(\mathbf{R}^{2n})$, if $m \in \mathbf{Z}$, $m > -n$;*

$$k_A(x,y) = G(x,y) \log|x-y| + F\Big(x, |x-y|, \frac{x-y}{|x-y|}\Big),$$

*where $F \in C^\infty(\mathbf{R}^n \times \overline{\mathbf{R}_+} \times \mathbf{S}^{n-1})$ and $G \in C^\infty(\mathbf{R}^{2n})$ vanishes to order $-n-m$ at $x = y$, if $m \in \mathbf{Z}$, $m \leq -n$; and*

$$k_A(x,y) = |x-y|^{-n-m} F\Big(x, |x-y|, \frac{x-y}{|x-y|}\Big) + G(x,y),$$

*where $F \in C^\infty(\mathbf{R}^n \times \overline{\mathbf{R}_+} \times \mathbf{S}^{n-1})$, $G \in C^\infty(\mathbf{R}^{2n})$, if $m \notin \mathbf{Z}$.*

*Proof.* This follows from the last proposition by a kind of argument which is very standard; we put it down in some detail here for the record but will be more brief on similar occasions later on.

In view of (14), the Schwartz kernel of $A$ is given by $k_A(x,y) = \check{a}(x, x-y)$, where $\check{a}(x,z) = \mathcal{F}^{-1}_{\xi \to z} a(x,\xi)$ is the inverse Fourier transform of the symbol $a$ of $A$ with respect to the $\xi$ variable. Using the polyhomogeneous expansion (13) of $a$, we can write, for any $N = 0, 1, 2, \ldots$,

$$a(x,\xi) = \theta(\xi) \sum_{j=0}^{N-1} a_{m-j}(x,\xi) + a_{(N)}(x,\xi),$$

where $\theta$ is a patch function, $a_{m-j} \in C^\infty(\mathbf{R}^n \times (\mathbf{R}^n \setminus \{0\}))$ is homogeneous in $\xi$ of degree $m-j$, and $a_{(N)} \in S^{m-N}$. Using the preceding proposition, this yields

$$\check{a}(x,z) = \chi(z) \sum_{j=0}^{N-1} \Big[U_{j-m-n}(x,z) + Q_{j-m-n}(x,z) \log|z|\Big] + g(x,z) + \check{a}_{(N)}(x,z),$$

with $\chi$, $U_{j-m-n}$, $Q_{j-m-n}$, $g$ as in the proposition, the latter three depending smoothly on $x$. By homogeneity, we can write

$$U_{j-m-n}(x,z) = |z|^{j-m-n} \mathfrak{u}_{j-m-n}(x, \tfrac{z}{|z|})$$

with $\mathfrak{u}_{j-m-n} \in C^\infty(\mathbf{R}^n \times \mathbf{S}^{n-1})$. A standard argument (imitating the proof of Borel's theorem, see e.g. [40], Proposition 3.5, or [29], Proposition 18.1.3) produces a function $\mathfrak{U}(x,r,\zeta)$, $\mathfrak{U} \in C^\infty(\mathbf{R}^n \times \mathbf{R} \times \mathbf{S}^{n-1})$, whose $j$-th order term in the Taylor expansion in the second variable is $r^j \mathfrak{u}(x,\zeta)$: that is,

$$\mathfrak{U}(x,r,\zeta) = \sum_{j=0}^{N-1} r^j \mathfrak{u}(x,\zeta) + O(r^N)$$

for any $N = 0, 1, 2, \ldots$. Similarly, one can find $\mathfrak{Q} \in C^\infty(\mathbf{R}^{2n})$ such that

$$\mathfrak{Q}(x,z) = \sum_{j=\max\{0,m+n\}}^{m+n+N-1} Q_{j-m-n}(x,z) + O(|z|^N)$$

for any $N = 0, 1, 2, \ldots$. Since $\check{a}_{(N)} \in C^{N-m-n-1}(\mathbf{R}^{2n})$ if $N \geq m+n+1$ (see e.g. [23], Theorem 8.8), and $N = 0, 1, 2, \ldots$ can be taken arbitrary, it follows that the function

$$c(x,z) = \check{a}(x,z) - |z|^{-m-n}\mathfrak{U}(x,|z|,\tfrac{z}{|z|}) - \mathfrak{Q}(x,z)\log|z|$$

belongs to $C^{N-m-n-1}(\mathbf{R}^{2n})$ for any $N$, hence, to $C^\infty(\mathbf{R}^{2n})$. Now $\mathfrak{Q} \equiv 0$ if $m \notin \mathbf{Z}$, while for $m \in \mathbf{Z}$ the function $c(x,z)$ can be combined together with the middle term on the right into $|z|^{-n-m}F(x,|z|,\tfrac{z}{|z|})$, and the assertion follows. $\square$

Note that since $x-y$ is a smooth function of $|x-y|$ and $\tfrac{x-y}{|x-y|}$, we get an equivalent statement if we replace $F(x,|x-y|,\tfrac{x-y}{|x-y|})$ in the last corollary by $F(y,|x-y|,\tfrac{x-y}{|x-y|})$, or even $F(x,y,|x-y|,\tfrac{x-y}{|x-y|})$ with $F \in C^\infty(\mathbf{R}^{2n} \times \overline{\mathbf{R}_+} \times \mathbf{S}^{n-1})$.

The last proof in fact establishes (omitting the variable $x$, which features only as a smooth parameter throughout) the validity of the first part of the following proposition; the second part is proved in the same way. Let $z = r\zeta$ ($r \in \overline{\mathbf{R}_+}, \zeta \in \mathbf{S}^{n-1}$) be the polar coordinates in $\mathbf{R}^n$. For $s \in \mathbf{R}$, let us denote by $\mathcal{C}_s$ the vector space of all distributions on $\mathbf{R}^n$ of the form

$$(r^s F(r,\zeta))\dot{} + G(r\zeta)\log r + P(D)\delta_0 + H,$$

where $F, G$ vanish for $r \geq 1$, $F \in C^\infty(\overline{\mathbf{R}_+} \times \mathbf{S}^{n-1})$, $H \in \mathcal{S}(\mathbf{R}^n)$, $P$ is a polynomial of degree $-s-n$ (thus $P = 0$ if $-s-n \notin \mathbf{N}$), and $G \in C^\infty(\mathbf{R}^n)$ vanishes to order $s$ at the origin if $s \in \mathbf{N}$, while $G = 0$ is $s \notin \mathbf{N}$. Here the prolongation $(r^s F)\dot{}$ of $r^s F$ from $\mathbf{R}^n \setminus \{0\}$ to $\mathbf{R}^n$ exists, since one can write the latter, choosing $k \in \mathbf{N}$ so large that $s + k > -n$, as the sum

$$r^s F(r,\zeta) = \sum_{j=0}^{k-1} \frac{r^{s+j}}{j!} \frac{\partial^j F(0,\zeta)}{\partial r^j} + r^{s+k} F^\#(r,\zeta)$$

of homogeneous functions of degrees $s, s+1, \ldots, s+k-1$ and a function of moderate growth in $L^1_{\text{loc}}(\mathbf{R}^n)$. Finally, let $\mathcal{P}_m$ denote the space of all polynomials of degree $\leq m$ on $\mathbf{R}^n$, and let again $\chi$ be a cutoff function and $\theta$ a patch function.

**Proposition 6.** (i) *The Fourier transform $\mathcal{F}$ is maps the Hörmander class $S^m$, $m \in \mathbf{R}$, into $\mathcal{C}_{-n-m}$.*

(ii) *In fact, $\mathcal{F}$ is a bijection of*

$$S^m + \theta \mathcal{P}_m \log := \{a + \theta p \log|\cdot| \,:\, a \in S^m, p \in \mathcal{P}_m\}$$

*onto $\mathcal{C}_{-n-m}$.*

*Proof.* (i) By definition, $S^m \subset \mathcal{S} + \theta \sum_{j \geq 0} \mathcal{H}_{m-j}$ (in the obvious sense, i.e. for $f \in S^m$ there exist $h_j \in \mathcal{H}_{m-j}$ such that $f - \theta \sum_{j=0}^{N-1} h_j \in S^{m-N}$ for all $N \in \mathbf{N}$;

in fact, the $h_j$ furthermore do not contain any logarithmic terms, nor the derivatives of $\delta_0$ which are killed by $\theta$). Using Proposition 4(i) and the same argument as in the last proof, this implies $\mathcal{F}S^m \subset \mathcal{S} + \chi \sum_{j \geq 0} \mathcal{H}_{-n-m+j} = \mathcal{C}_{-n-m}$.

(ii) Now $S^m + \theta \mathcal{P}_m \log = \mathcal{S} + \theta \sum_{j \geq 0} \mathcal{H}_{m-j}$ (this time we use also the logarithmic terms in $\mathcal{H}_{m-j}$, while the derivatives of $\delta_0$ still play no role); since $\mathcal{C}_{-n-m} = \mathcal{S} + \chi \sum_{j \geq 0} \mathcal{H}_{-n-m+j}$, the claim follows in the same way using part (ii) of Proposition 4. □

Using local charts, one gets also the analogue of Corollary 5 for operators on $\partial \Omega$.

**Corollary 7.** *If $\Omega \subset \mathbf{R}^n$ is a bounded domain with smooth boundary and $A \in \Psi^m_{\mathrm{cl}}(\partial \Omega)$, then its Schwartz kernel (with respect to some smooth volume element on $\partial \Omega$) is of the same form as in Corollary 5, with $F \in C^\infty(\partial \Omega \times \overline{\mathbf{R}_+} \times \mathbf{S}^{n-1})$.*

*Proof.* Immediate from Corollary 5 and the observation that, for any diffeomorphism $\Phi$, $\frac{|\Phi(x)-\Phi(y)|}{|x-y|}$ and $\frac{\Phi(x)-\Phi(y)}{|\Phi(x)-\Phi(y)|}$ are smooth functions of $x$, $|x-y|$ and $\frac{x-y}{|x-y|}$. □

Let us now give the same treatment to boundary $\Psi$DOs.

**4.2 Potential operators.** We start with the Poisson type operators.

**Theorem 8.** *If $K \in \mathcal{K}^d(\overline{\Omega})$, where $\Omega \subset \mathbf{R}^n$ is a bounded domain with smooth boundary, or $\Omega = \mathbf{R}^n_+$, then its Schwartz kernel $k_K$ satisfies*

$$k_K(x,\zeta) = |x-\zeta|^{1-n-d} F\Big(\zeta, |x-\zeta|, \frac{x-\zeta}{|x-\zeta|}\Big) + G(x,\zeta) \log|x-\zeta|$$

*if $d \in \mathbf{Z}$, $d > 1-n$;*

$$k_K(x,\zeta) = G(x,\zeta) \log|x-\zeta| + F\Big(\zeta, |x-\zeta|, \frac{x-\zeta}{|x-\zeta|}\Big)$$

*with $G(x,\zeta)$ vanishing to order $1-n-d$ on $x=\zeta$, if $d \in \mathbf{Z}$, $d \leq 1-n$; and*

$$k_K(x,\zeta) = |x-\zeta|^{1-n-d} F\Big(\zeta, |x-\zeta|, \frac{x-\zeta}{|x-\zeta|}\Big) + G(x,\zeta)$$

*if $d \notin \mathbf{Z}$, for $x \in \Omega$, $\zeta \in \partial \Omega$, where $F \in C^\infty(\partial \Omega \times \overline{\mathbf{R}_+} \times \mathbf{S}^{n-1})$, $G \in C^\infty(\overline{\Omega} \times \partial \Omega)$.*

Note that by Seeley's extension theorem, we can even assume that $F \in C^\infty(\mathbf{R}^n \times \mathbf{R} \times \mathbf{R}^n)$ and $G \in C^\infty(\mathbf{R}^n \times \mathbf{R}^n)$.

*Proof.* It is enough to deal with the case of $\Omega = \mathbf{R}^n_+$ and $K$ of the form (17) (the negligible piece of the operator contributes a Schwartz kernel smooth on $\overline{\mathbf{R}^n_+} \times \overline{\mathbf{R}^n_+}$), where it is now convenient to take the $y'$-version of the kernel $k$, i.e.

$$Ku(x) = \int_{\mathbf{R}^{n-1}} \int_{\mathbf{R}^{n-1}} e^{i(x'-y') \cdot \xi'} k(y', x_n, \xi') u(y') \, dy' d\!\!\!{}^-\!\xi'$$

interpreted as oscillatory integral, with $k \in S^{d-1}_{\mathrm{cl}}(\mathbf{R}^{n-1}, \mathbf{R}^{n-1}, \mathcal{S}_+)$. The Schwartz kernel is thus given by

$$k_K(x, y') = \check{k}(y', x_n, x'-y'),$$

where $\check{k}$ denotes the inverse Fourier transform of $k(y', x_n, \xi')$ with respect to $\xi'$.

We claim, first of all, that

$$(36) \qquad k \in S_{\mathrm{cl}}^{d-1} \implies \check{k} \in C^j(\mathbf{R}^{n-1} \times \overline{\mathbf{R}_+^n}) \qquad \text{if } d < 1 - n - j.$$

Indeed, by the definition of $S_{\mathrm{cl}}^{d-1}$, the hypothesis implies that for any $l \in \mathbf{N}$, $\beta \in \mathbf{N}^{n-1}$ and compact subset $K$ of $\mathbf{R}^{n-1}$,

$$\sup_{y' \in K, x_n \geq 0, \xi' \in \mathbf{R}^{n-1}} \frac{|\partial_{x_n}^l \xi'^\beta k(y', x_n, \xi')|}{(1 + |\xi'|)^{d+l+|\beta|}} < \infty.$$

Thus by the Riemann-Lebesgue lemma, $\partial_{x_n}^l \partial_{x'}^\beta \check{k}(y', x_n, x')$ is continuous and uniformly bounded as $y' \in K$, $x \in \overline{\mathbf{R}_+^n}$ if $d + l + |\beta| < 1 - n$, proving the claim.

Next, recalling the polyhomogeneous expansion

$$(37) \qquad k = \sum_{j=0}^{N-1} k_{d-1-j} + k^{(N)}, \qquad k^{(N)} \in S_{\mathrm{cl}}^{d-1-N}, \qquad N = 0, 1, 2, \ldots,$$

of $k$, we claim that

$$(38) \qquad \check{k}_{d-1-j}(y', x) = U_{1-n-d+j}(y', x) + Q_{1-n-d+j}(y', x) \log |x| + s_j(y', x),$$

with $U_{1-n-d+j} \in C^\infty(\mathbf{R}^{n-1} \times \overline{\mathbf{R}_+^n})$ homogeneous in $x$ and $Q_{1-n-d+j} \in C^\infty(\mathbf{R}^{n-1} \times \mathbf{R}^n)$ a homogeneous polynomial in $x$ of the degree indicated by the subscript, and $s_j \in C^\infty(\mathbf{R}^{n-1} \times \overline{\mathbf{R}_+^n})$. Combining this with (37) and (36) leads, exactly as in the proof of Corollary 5, to the formulas in the theorem for $\Omega = \mathbf{R}_+^n$, while passage to local charts as in the proof of Corollary 7 finally yields the result in general.

From now on, we thus assume (replacing $d - 1 - j$ by $d - 1$ for brevity) that $k = k_{d-1}$ satisfies the homogeneity condition

$$(39) \qquad k(y', \frac{x_n}{t}, t\xi') = t^d k(y', x_n, \xi') \qquad \text{for } \lambda \geq 1, |\xi'| \geq 1.$$

Since $y'$ enters — from the point of view of our inverse Fourier transform in $\xi'$ — only as a parameter on which everything depends smoothly, we will drop it from the notation for the rest of the proof. Define $l(x_n, \xi') := |\xi'|^d k(|\xi'| x_n, \xi'/|\xi'|)$, or, in the polar coordinates $\xi' = r\zeta$ ($r \geq 0, \zeta \in \mathbf{S}^{n-2}$) on $\mathbf{R}^{n-1}$

$$l(x_n, r\zeta) := r^d k(r x_n, \zeta).$$

Then $l \in C^\infty(\overline{\mathbf{R}_+} \times (\mathbf{R}^{n-1} \setminus \{0\}))$,

$$(40) \qquad \begin{aligned} l(x_n, \xi') &= k(x_n, \xi') \qquad \text{if } |\xi'| \geq 1, \\ l(\tfrac{x_n}{t}, t\xi') &= t^d l(x_n, \xi') \qquad \forall t > 0, \xi' \neq 0, x_n \geq 0. \end{aligned}$$

For brevity, we also denote $l_t(\xi') := l(t, \xi')$; introducing the dilation operator

$$\delta_t f(x) := f(tx), \qquad t > 0,$$

we thus have
$$l_t = t^{-d}\delta_t l_1.$$

Since $k \in C^\infty(\overline{\mathbf{R}^n_+})$ belongs to $\mathcal{S}_+$ as a function of $x_n$ for each fixed $\xi'$, it follows that $\theta l_1 \in \mathcal{S}_+$ for any patch function $\theta$, while

$$(41) \qquad l_1(r\zeta) = \sum_{j=0}^{N-1} \frac{r^{d+j}}{j!} \Big[\frac{\partial^j k(x_n,\zeta)}{\partial x_n^j}\Big]_{x_n=0} + O(r^{d+N}) \qquad \text{as } r \searrow 0.$$

Thus in the notation of §4.1, $l_1 \in \mathcal{C}_d(\mathbf{R}^{n-1})$.

Assume first that $d > 1 - n$, so that $l_1$ — as well as any $l_t$, $t \geq 0$ — is integrable as the origin, and thus defines a distribution on all of $\mathbf{R}^{n-1}$. For its inverse Fourier transform, we have by Proposition 6(ii)

$$\check{l}_1 \in \mathcal{C}_d^\vee = S^{1-n-d} + \theta\mathcal{P}_{1-n-d}\log = S^{1-n-d},$$

and

$$\check{l}_{x_n}(x) = (x_n^{-d}\delta_{x_n}l_1)^\vee(x) = x_n^{1-n-d}\check{l}_1\Big(\frac{x'}{x_n}\Big), \qquad x_n > 0.$$

The right-hand side of the last formula is a smooth function on $\mathbf{R}^n_+$ homogeneous of degree $1 - n - d$, and from $\check{l}_1 \in S^{1-n-d}$ it follows that it actually extends smoothly to $\overline{\mathbf{R}^n_+} \setminus \{0\}$: namely, if $\check{l}^1$ has the polyhomogeneous expansion

$$\check{l}_1(x') \sim \sum_{j=0}^\infty |x'|^{1-n-d-j} q_j(\tfrac{x'}{|x'|}),$$

then

$$\check{l}_{x_n}(x') \sim \sum_{j=0}^\infty x_n^j |x'|^{1-n-d-j} q_j(\tfrac{x'}{|x'|}).$$

We now claim that

$$(42) \qquad \check{k}(x_n, x') - \check{l}_{x_n}(x') =: s_0(x) \quad \text{satisfies } s_0 \in C^\infty(\overline{\mathbf{R}^n_+}).$$

Setting $U_{1-n-d}(x_n, x') = \check{l}_{x_n}(x')$ we then get (38) (with $Q_{1-n-d} \equiv 0$), thus settling the case of $d > 1 - n$.

It remains to show (42). Fix a cutoff function $\chi$ with $\chi(\xi') = 1$ for $|\xi'| \leq 1$ and $\chi(\xi') = 0$ for $|\xi'| \geq 2$, and let $\theta = 1 - \chi$ be the corresponding patch function. Then $\theta k = \theta l$, so

$$s_0 = (\chi k - \chi l)^\vee.$$

Differentiating under the integral sign in

$$(\chi k)^\vee(x_n, x') = \int_{|\xi'|<2} e^{ix'\cdot\xi'}\chi(\xi')k(x_n,\xi')đ\xi'$$

and noting that $k \in C^\infty(\overline{\mathbf{R}_+} \times \mathbf{R}^{n-1})$ and $k(\cdot,\xi') \in \mathcal{S}_+$ for each fixed $\xi'$, hence $\partial_{x_n}^m k(x_n,\xi') \in L^\infty(\overline{\mathbf{R}_+} \times K)$ for any compact $K \subset \mathbf{R}^{n-1}$, it is immediate that

$$(\chi k)^\vee \in C^\infty(\overline{\mathbf{R}^n_+}).$$

On the other hand, from $l(x_n, r\zeta) = r^d k(rx_n, \zeta)$ and the fact (again) that $k \in C^\infty(\overline{\mathbf{R}_+} \times \mathbf{R}^{n-1})$ and $k(\cdot, \zeta) \in \mathcal{S}_+$ for each fixed $\zeta \in \mathbf{S}^{n-2}$, we see that

$$\partial_{x_n}^m l(x_n, r\zeta) = r^{d+m}(\partial_{x_n}^m k)(rx_n, \zeta)$$

with $\partial_{x_n}^m k$ bounded on $\overline{\mathbf{R}_+} \times \mathbf{S}^{n-2}$. Differentiating under the integral sign in

$$(\chi l)^\vee(x_n, x') = \int_{|\xi'|<2} e^{ix'\cdot\xi'} \chi(\xi') l(x_n, \xi') d\!\!\!/\xi',$$

it again transpires that

$$(\chi l)^\vee \in C^\infty(\overline{\mathbf{R}_+^n}).$$

This completes the proof of (42) and, consequently, (38) for the case of $d > 1 - n$.

For $d \leq 1 - n$, one can reduce to the previous case upon differentiating $m$ times with respect to $x_n$, where $m = [2-n-d]$: the kernel $\partial_{x_n}^m k$ belongs to $S_{\text{cl}}^{d-1-m}$ and satisfies (39) with $d$ replaced by $d + m > 1 - n$. By what we have already proved, we thus have

$$\partial_{x_n}^m \check{k} = (\partial_{x_n}^m k)^\vee = U_{1-n-d-m} + s_0, \tag{43}$$

with $U \in C^\infty(\overline{\mathbf{R}_+^n})$ homogeneous of the indicated degree, and $s_0 \in C^\infty(\mathbf{R}^n)$.

Observe now that if $f \in C^\infty(\overline{\mathbf{R}_+^n} \setminus \{0\})$ is homogeneous of degree $s < -1$, then $\int_\infty f(x)\,dx_n$ (exists and) is homogeneous of degree $s + 1$; if $f$ is homogeneous of degree $s > -1$, then $\int_0 f(x)\,dx_n$ is homogeneous of degree $s + 1$; and if $f$ is homogeneous of degree $-1$, then expressing the difference $f(x) - ax_n|x|^{-2}$, where $a = f(0,1)$, on the hyperplane $x_n = 1$ as

$$f(x', 1) - \frac{a}{1+|x'|^2} = \chi(x') \sum_{j=1}^{n-1} x_j f_j(x') + |x'|^2 g(x')$$

with $f_j, g \in C^\infty(\mathbf{R}^{n-1})$ and $\chi$ a cutoff function on $\mathbf{R}^{n-1}$ (this is achieved simply by taking any $f_j$ such that $\sum_{j=1}^{n-1} x_j f_j(x')$ has the same Taylor expansion around $x' = 0$ as the left-hand side), we get

$$f(x) = \frac{ax_n}{|x|^2} + \sum_{j=1}^{n-1} x_j F_j(x) + |x'|^2 G(x), \tag{44}$$

where the functions

$$F_j(x) := \frac{1}{x_n^2} \chi\!\left(\frac{x'}{x_n}\right) f_j\!\left(\frac{x'}{x_n}\right), \qquad G(x) := \frac{1}{x_n^3} g\!\left(\frac{x'}{x_n}\right)$$

are smooth on $\overline{\mathbf{R}_+^n} \setminus \{0\}$ (this is clear for $F_j$, and follows from (44) for $G$) and homogeneous there of degrees $-2$ and $-3$, respectively, so the function

$$F(x) := a \log|x| + \sum_{j=1}^{n-1} x_j \int_\infty F_j(x)\,dx_n + |x'|^2 \int_\infty G(x)\,dx_n$$

belongs to $C^\infty(\overline{\mathbf{R}^n_+} \setminus \{0\})$, $\frac{\partial F}{\partial x_n} = f$ and $F(x) - a \log |x|$ is homogeneous of degree 0.

Integrating (43) with respect to $x_n$ in the above way $m$ times, it transpires that $\check{k}$ has the form

$$\text{(45)} \quad U_{1-n-d}(x) + Q_{1-n-d}(x) \log |x| + S_0(x) + \sum_{j=0}^{m-1} x_n^j v_j(x'),$$

where $U_{1-n-d} \in C^\infty(\overline{\mathbf{R}^n_+})$ is homogeneous of the indicated degree, $Q_{1-n-d}$ is a polynomial of the indicated degree (in particular, $Q \equiv 0$ if $1 - n - d \notin \mathbf{N}$), $S_0 := (\int_0 dx_n)^m s_0 \in C^\infty(\mathbf{R}^n)$, and $v_j$ are some distributions on $\mathbf{R}^{n-1}$. Since $\check{k}$ is clearly a Schwartz function for each fixed $x_n > 0$, so must be (45), implying that in fact $v_j \in C^\infty(\mathbf{R}^{n-1})$. Setting $s_m := S_0 + \sum_{j=0}^{m-1} x_n^j v_j \in C^\infty(\mathbf{R}^n)$, we thus get (38). This settles the case $d \leq 1 - n$, and thus completes the proof of the whole theorem. $\square$

*Remark 9.* Note that the last proof gives actually a somewhat more precise information than (38): namely, if $k_{d-1} \in S_{\text{cl}}^{d-1}(\mathbf{R}^{n-1}, \mathbf{R}^{n-1}, \mathcal{S}_+)$ is of the form

$$k_{d-1} = x_n^p k_{d+p-1} \quad \text{with some } p \in \mathbf{N} \text{ and } k_{d+p-1} \in S_{\text{cl}}^{d+p-1}(\mathbf{R}^{n-1}, \mathbf{R}^{n-1}, \mathcal{S}_+),$$

then

$$\check{k}_{d-1}(y', x) = x_n^p U_{1-n-d-p}(y', x) + x_n^p Q_{1-n-d-p}(y', x) \log |x| + x_n^p s(y', x)$$

with $U_{1-n-d-p} \in C^\infty(\mathbf{R}^{n-1} \times \overline{\mathbf{R}^n_+})$ homogeneous in $x$ and $Q_{1-n-d-p} \in C^\infty(\mathbf{R}^{n-1} \times \mathbf{R}^n)$ a homogeneous polynomial in $x$ of the indicated degree, and $s \in C^\infty(\mathbf{R}^{n-1} \times \overline{\mathbf{R}^n_+})$. That is, loosely speaking, if it is possible to pull out a factor of $x_n^p$ from $k$, then the same is true for $\check{k}$. This is immediate from the obvious fact that the inverse Fourier transform with respect to $\xi'$ commutes with multiplication by $x_n$. It will prove useful in Section 5. $\square$

**4.3 Singular Green operators.** The Schwartz kernel of trace operators of class zero can be obtained simply by taking the adjoint of potential operators; the case of general class $r \geq 0$ then follows easily from (18). It therefore remains to deal with singular Green operators.

Denote by $\overline{\mathbf{S}^n_+} = \mathbf{S}^n \cap \overline{\mathbf{R}^{n+1}_+}$ the closed upper hemisphere in $\mathbf{R}^{n+1}$.

**Theorem 10.** *Let $A \in \mathcal{G}_0^d(\overline{\mathbf{R}^n_+})$ and denote by $\tilde{y} = (y', -y_n)$ the image of $y = (y', y_n) \in \mathbf{R}^n$ under the reflection with respect to the hyperplane $y_n = 0$. Then the Schwartz kernel $k_A$ of $A$ satisfies*

$$k_A(x, y) = |z|^{-n-d} F\Big(x', |z|, \frac{z}{|z|}\Big) + G(x', z) \log |z|$$

*if $d \in \mathbf{Z}$, $d > -n$;*

$$k_A(x, y) = F\Big(x', |z|, \frac{z}{|z|}\Big) + G(x', z) \log |z|,$$

*with $G$ vanishing to order $-n - d$ at $z = 0$, if $d \in \mathbf{Z}$, $d \leq -n$; and*

$$k_A(x, y) = |z|^{-n-d} F\Big(x', |z|, \frac{z}{|z|}\Big) + G(x', z)$$

if $d \notin \mathbf{Z}$, for $x, y \in \mathbf{R}_+^n$, where

$$z = z(x, y) := (x' - y', x_n, y_n) \in \mathbf{R}^{n+1}$$

and $F \in C^\infty(\mathbf{R}^{n-1} \times \overline{\mathbf{R}_+} \times \overline{\mathbf{S}_+^n})$, $G \in C^\infty(\mathbf{R}^{n-1} \times \mathbf{R}^{n+1})$.

In particular, $k_A$ extends smoothly up to the boundary of $\overline{\mathbf{R}_+^n} \times \overline{\mathbf{R}_+^n}$ away from the boundary diagonal $x = \tilde{y}$.

Notice that an equivalent statement is again obtained upon replacing $F(x', |z|, \frac{z}{|z|})$ by $F(x', y', |z|, \frac{z}{|z|})$, or even $F(x, y, |z|, \frac{z}{|z|})$, and similarly for $G$. Also, by Seeley's extension theorem, one may replace $\overline{\mathbf{S}_+^n}$ by $\mathbf{S}^n$ or even $\mathbf{R}^{n+1}$.

*Proof.* The idea of the proof closely parallels that for Theorem 8, so we will be more brief. By (20), where we again take the kernel $g$ in the $y'$-form rather than the $x'$-form,

$$k_A(x, y) = \check{g}(y', x_n, y_n, x' - y')$$

where $\check{g}$ denotes the inverse Fourier transform of $g(y', x_n, y_n, \xi')$ with respect to $\xi'$. An argument involving the Riemann-Lebesgue lemma again shows that $\check{g} \in C^j(\overline{\mathbf{R}_+^n} \times \overline{\mathbf{R}_+^n})$ if $d < -n - j$, and as before it follows that it is enough to deal with each term of the polyhomogeneous expansion of $g$ separately. We thus assume from now on that $g \in S_{cl}^{d-1}(\mathbf{R}^{n-1}, \mathbf{R}^{n-1}, \mathcal{S}_{++})$ satisfies

(46) $\qquad g(y', \frac{x_n}{t}, \frac{y_n}{t}, t\xi') = t^{d+1} g(y', x_n, y_n, \xi') \qquad \text{for } t \geq 1, |\xi'| \geq 1.$

Omitting again $y'$ (which plays only the role of a smooth parameter throughout) from the notation, set for $\xi' = r\zeta$, $r \in \overline{\mathbf{R}_+}$, $\zeta \in \mathbf{S}^{n-1}$,

$$l(x_n, y_n, r\zeta) := r^{d+1} g(rx_n, ry_n, \zeta).$$

Then $l \in C^\infty(\overline{\mathbf{R}_+} \times \overline{\mathbf{R}_+} \times (\mathbf{R}^{n-1} \setminus \{0\}))$,

$$l(x_n, y_n, \xi') = g(x_n, y_n, \xi') \qquad \text{if } |\xi'| \geq 1,$$

and

(47) $\qquad l(\frac{x_n}{t}, \frac{y_n}{t}, t\xi') = t^{d+1} l(x_n, y_n, \xi') \qquad \forall t > 0, \xi' \neq 0.$

We further denote

$$l_{R,\vartheta}(\xi') := l(R \cos \vartheta, R \sin \vartheta, \xi'), \qquad R \in \overline{\mathbf{R}_+}, \vartheta \in [0, \frac{\pi}{2}], \xi' \in \mathbf{R}^{n-1} \setminus \{0\};$$

in terms of the dilation operator $\delta_t f(x) := f(tx)$, we thus have

$$l_{R,\vartheta} = R^{-d-1} \delta_R l_{1,\vartheta}.$$

Since $g \in C^\infty(\overline{\mathbf{R}_+} \times \overline{\mathbf{R}_+} \times \mathbf{R}^{n-1})$ belongs to $\mathcal{S}_{++}$ as a function of $x_n, y_n$ for each fixed $\xi'$, which implies that $g(R \cos \vartheta, R \sin \vartheta, \xi')$ belongs to $\mathcal{S}_+$ as a function of $R$ for each fixed $\xi'$ and $\vartheta$, it follows that $\theta l_{1,\vartheta} \in \mathcal{S}_+$ for any patch function $\theta$, while

$$l_{1,\vartheta}(r\zeta) = \sum_{j+k \leq N-1} \frac{r^{d+1+j+k}}{j! k!} (\cos \vartheta)^j (\sin \vartheta)^k [\partial_{x_n}^j \partial_{y_n}^k g](0, 0, \zeta) + O(r^{N+d+1})$$

as $r \searrow 0$. Thus in the notation of §4.1, $l_{1,\vartheta} \in \mathcal{C}_{d+1}(\mathbf{R}^{n-1})$.

Assume first that $d > -n$, so that $l_{1,\vartheta}$ — as well as any $l_{R,\vartheta}$, $R \geq 0$ — is integrable at the origin, and thus defines a distribution on all of $\mathbf{R}^{n-1}$. For its inverse Fourier transform, we have by Proposition 6(ii)

$$\check{l}_{1,\vartheta} \in \mathcal{C}_{d+1}^{\vee} = S^{-n-d} + \theta \mathcal{P}_{-n-d} \log = S^{-n-d}$$

(since $-n - d < 0$), and

$$\check{l}_{R,\vartheta}(x') = (R^{-d-1}\delta_R l_{1,\vartheta})^{\vee}(x') = R^{-n-d}\check{l}_{1,\vartheta}\Big(\frac{x'}{R}\Big), \qquad R > 0.$$

Since $l_{1,\vartheta}$ depends smoothly on $\vartheta$, the right-hand side of the last formula is a smooth function of $(x_n, y_n, x') = (R \cos \vartheta, R \sin \vartheta, x') \in \mathbf{R}_+ \times \mathbf{R}_+ \times \mathbf{R}^{n-1}$, homogeneous of degree $-n - d$, and from $\check{l}_{1,\vartheta} \in S^{-n-d}(\mathbf{R}^{n-1})$ it follows that it actually extends smoothly to $\overline{\mathbf{R}_+} \times \overline{\mathbf{R}_+} \times \mathbf{R}^{n-1}$: namely, if $\check{l}_{1,\vartheta}$ has the polyhomogeneous expansion

$$\check{l}_{1,\vartheta}(x') \sim \sum_{j=0}^{\infty} |x'|^{-n-d-j} q_j\Big(\frac{x'}{|x'|}, \vartheta\Big),$$

then

$$\check{l}_{R,\vartheta}(x') \sim \sum_{j=0}^{\infty} R^j |x'|^{-n-d-j} q_j\Big(\frac{x'}{|x'|}, \vartheta\Big),$$

so $\check{l}_{R,\vartheta}(x')$ is $C^\infty$ in $(R, \vartheta, x') \in \overline{\mathbf{R}_+} \times [0, \frac{\pi}{2}] \times \mathbf{R}^{n-1}$. Arguing as in the proof of Theorem 8, one finally checks that for any patch function $\theta$ on $\mathbf{R}^{n-1}$ vanishing on $|\xi'| \leq 1$, and $\chi = 1 - \theta$ the corresponding cutoff function, both $(\chi g)^{\vee}$ and $(\chi l)^{\vee}$ and, hence, also $(\chi g - \chi l)^{\vee} = (k - l)^{\vee}$, are smooth on $\overline{\mathbf{R}_+} \times \overline{\mathbf{R}_+} \times \mathbf{R}^{n-1}$. Thus, setting $U_{-n-d}(y_n, x) = \check{l}(x_n, y_n, x')$ and restoring the variable $y'$, we have proved that for $d > -n$ and $g \in S_{\mathrm{cl}}^{d-1}$ satisfying (46),

$$\check{g}(y', x_n, y_n, x') = U_{-n-d}(y', z) + s(y', z),$$
$$z = (x_n, y_n, x') \in \overline{\mathbf{R}_+} \times \overline{\mathbf{R}_+} \times \mathbf{R}^{n-1} =: \overline{\mathbf{R}_{++}^{n+1}}, y' \in \mathbf{R}^{n-1},$$

with $U_{-n-d} \in C^\infty(\mathbf{R}^{n-1} \times (\overline{\mathbf{R}_{++}^{n+1}} \setminus \{0\}))$ homogeneous of the indicated degree in $z$, and $s \in C^\infty(\mathbf{R}^{n-1} \times \overline{\mathbf{R}_{++}^{n+1}})$.

Finally, the case of $d \leq -n$ is again handled by differentiating $g$ with respect to $x_n$ (or $y_n$) $m$ times, where $m = [1 - n - d]$, and then integrating $m$ times exactly as in the proof of Theorem 8. The conclusion is that, for any $d \in \mathbf{R}$,

$$\check{g}(y', x_n, y_n, x') = U_{-n-d}(y', z) + Q_{-n-d}(y', z) \log |z| + s(y', z),$$
(48)
$$z = (x_n, y_n, x') \in \overline{\mathbf{R}_{++}^{n+1}}, y' \in \mathbf{R}^{n-1},$$

where $U \in C^\infty(\mathbf{R}^{n-1} \times (\overline{\mathbf{R}_{++}^{n+1}} \setminus \{0\}))$ is homogeneous in $z$ and $Q_{-n-d} \in C^\infty(\mathbf{R}^{n-1} \times \mathbf{R}^{n+1})$ is a homogeneous polynomial in $z$ of the indicated degrees, and $s \in C^\infty(\mathbf{R}^{n-1} \times \overline{\mathbf{R}_{++}^{n+1}})$. Putting everything together, the assertion follows in the same way as in the proof of Corollary 5. $\square$

*Remark* 11. Just as in Remark 9, the last proof gives a bit more precise information: namely, if a factor of $x_n^p y_n^q$, with some $p, q \in \mathbf{N}$, can be pulled out from $g \in S_{\text{cl}}^{d-1}$ in the sense that

$$g = x_n^p y_n^q G \qquad \text{with } G \in S_{\text{cl}}^{d+p+q-1}(\mathbf{R}^{n-1}, \mathbf{R}^{n-1}, \mathcal{S}_{++}),$$

then it can also be pulled out from the $U$, $Q$ and $s$ in (48). $\square$

Note that in the formulas in Theorem 10, one can replace $|z|$ and $F(x', |z|, \frac{z}{|z|}) \in C^\infty(\mathbf{R}^{n-1}, \overline{\mathbf{R}_+}, \overline{\mathbf{S}_+^n})$ by $|x - \tilde{y}|$ and $F(x, |x - \tilde{y}|, \frac{x-\tilde{y}}{|x-\tilde{y}|}) \in C^\infty(\overline{\mathbf{R}_+^n} \times \overline{\mathbf{R}_+} \times \mathbf{S}^n)$, respectively, since the ratio

$$(49) \qquad \frac{|x - \tilde{y}|^2}{|z|^2} = \frac{|x' - y'|^2 + (x_n + y_n)^2}{|x' - y'|^2 + x_n^2 + y_n^2}$$

extends to a smooth positive function of $(x, |x-\tilde{y}|, \frac{x-\tilde{y}}{|x-\tilde{y}|}) \in \overline{\mathbf{R}_+^n} \times \overline{\mathbf{R}_+} \times \mathbf{S}^n$. Similarly, $|z|$ can be replaced by $d_\Delta(x, y)$, the distance of $(x, y) \in \overline{\mathbf{R}_+^n} \times \overline{\mathbf{R}_+^n}$ to the boundary diagonal $\{(x, x) : x_n = 0\}$.

The last theorem is again easily transferred from $\overline{\mathbf{R}_+^n}$ to manifolds by means of local charts. Recall that for $\Omega$ a bounded domain in $\mathbf{R}^n$ with smooth boundary, for $\epsilon > 0$ small enough the mapping

$$\pi : \partial\Omega \times (-\epsilon, +\epsilon) \to \mathbf{R}^n, \qquad \pi(\zeta, t) = \zeta + t n_\zeta,$$

where $n_\zeta$ is the unit inward normal vector at $\zeta \in \partial\Omega$, is a diffeomorphism. Denote its image by $\mathcal{V}_\epsilon$, and define the map $y \mapsto \tilde{y} \colon \mathcal{V}_\epsilon \to \mathcal{V}_\epsilon$ by

$$\widetilde{\pi(\zeta, t)} = \pi(\zeta, -t)$$

(the "reflection with respect to $\partial\Omega$").

**Corollary 12.** *Let $\Omega$ be a bounded domain in $\mathbf{R}^n$ with smooth boundary, $\epsilon$ be as above, and $A \in \mathcal{G}_0^d(\overline{\Omega})$. Then the Schwartz kernel $k_A$ of $A$ is $C^\infty$ outside $\mathcal{V}_{\epsilon/2} \times \mathcal{V}_{\epsilon/2}$, while on $\mathcal{V}_\epsilon \times \mathcal{V}_\epsilon$ it satisfies*

$$k_A(x, y) = |w|^{-n-d} F\Big(x, y, |w|, \frac{w}{|w|}\Big) + G(x, y) \log |w|$$

*if $d \in \mathbf{Z}$, $d > -n$;*

$$k_A(x, y) = F\Big(x, y, |w|, \frac{w}{|w|}\Big) + G(x, y) \log |w|,$$

*with $G$ vanishing to order $-n - d$ at $x = \tilde{y}$, if $d \in \mathbf{Z}$, $d \leq -n$; and*

$$k_A(x, y) = |w|^{-n-d} F\Big(x, y, |w|, \frac{w}{|w|}\Big) + G(x, y)$$

*if $d \notin \mathbf{Z}$, where*

$$w = x - \tilde{y}$$

*and $F \in C^\infty(\overline{\Omega} \times \overline{\Omega} \times \overline{\mathbf{R}_+} \times \mathbf{S}^{n-1})$, $G \in C^\infty(\overline{\Omega} \times \overline{\Omega})$.*

Again, $F(x,y,|w|,\frac{w}{|w|})$ could be replaced by $F(x,|w|,\frac{w}{|w|})$ or $F(y,|w|,\frac{w}{|w|})$, or extended smoothly to all of $\mathbf{R}^{n+n+1+n}$.

*Proof.* Note first of all that all this holds in the local chart $\overline{\mathbf{R}^n_+}$, by the last theorem (and the remarks after it), since $z = (x', x_n, y_n)$ is a smooth function of $x, y$, while $|z|$ is a smooth function of $x, y, |w|$ and $\frac{w}{|w|}$ (cf. (49)). The general case now follows as in the proof of Corollary 7. □

Finally, the Schwarz kernel of an arbitrary singular Green operator of class $r \geq 0$ can be written down easily using the last corollary, (19) and Theorem 8; we omit the details.

## 5. The Poisson kernel

Recall that the kernel $k(x', x_n, \xi') \in S^{d-1}_{\text{cl}}(\mathbf{R}^{n-1}, \mathbf{R}^{n-1}, \mathcal{S}_+)$ of a potential operator (17) has the polyhomogeneous expansion

$$k(x', x_n, \xi') \sim \sum_{j=0}^{\infty} k_{d-j}(x', x_n, \xi') \tag{50}$$

where $k_{d-j} \in S^{d-j-1}_{\text{cl}}(\mathbf{R}^{n-1}, \mathbf{R}^{n-1}, \mathcal{S}_+)$ possess the homogeneity property

$$k_{d-j}(x', \tfrac{x_n}{\lambda}, \lambda\xi') = \lambda^{d-j} k_{d-j}(x', x_n, \xi') \tag{51}$$

for $\lambda \geq 1$ and $|\xi'| \geq 1$. It is sometimes convenient — and we will do that from now on — to redefine the $k_{d-j}$ by homogeneity on $|\xi'| < 1$ so that (51) actually holds for all $\lambda > 0$ and $\xi' \neq 0$; of course, $k_{d-j}$ will then in general have a singularity as $\xi' = 0$ (i.e. will belong only to $C^\infty(\mathbf{R}^{n-1} \times \overline{\mathbf{R}_+} \times (\mathbf{R}^{n-1} \setminus \{0\})))$, and also the expansion (50) will then hold only for $|\xi'| \geq 1$, i.e.

$$\sup_{x \in K,\, x_n \geq 0,\, |\xi'| \geq 1} \frac{\left| x_n^k \partial_{x_n}^l \partial_{x'}^\beta \partial_{\xi'}^\alpha \left[ k - \sum_{j=0}^{N-1} k_{d-j} \right] \right|}{|\xi|^{d-N+l-k-|\alpha|}} < \infty \tag{52}$$

for all $k, l, N \geq 0$, multiindices $\alpha, \beta$ and compact subsets $K$ of $\mathbf{R}^{n-1}$.

We have seen in the proof of Theorem 8 in §4.2 above that each homogeneous component

$$U_{1-n-d+j}(x', z) + Q_{1-n-d+j}(x', z) \log|z| \tag{53}$$

in the Schwartz kernel $k_K(x, y') = \check{k}(x', z)$, $z = (x_n, x' - y')$, of the potential operator in question (with $U_{1-n-d+j}, Q_{1-n-d+j}$ homogeneous of the indicated degree in $z$ — the latter a polynomial — and smoothly depending on $x'$) actually arises, modulo an error term smooth on all of $\mathbf{R}^{n-1} \times \overline{\mathbf{R}^n_+}$, precisely as the inverse Fourier transform with respect to $\xi'$ of

$$k_{d-j}(x', x_n, \xi')\theta(\xi')$$

with some patch function $\theta$. We use this correspondence to obtain some extra information about the functions $U_{1-n-d+j}, Q_{1-n-d+j}$ in (53) in the particular case of our Poisson operator $\mathbf{K}$.

So let us consider our Dirichlet problem

$$\Delta f = 0, \qquad f|_{\partial\Omega} = u$$

on a bounded domain $\Omega \subset \mathbf{R}^n$ with smooth boundary $\partial\Omega$. We want to compute the components (53) of the corresponding Poisson kernel $K(x,\zeta)$ for $x,\zeta$ near some boundary point $a \in \partial\Omega$. Without loss of generality, we may assume that $a$ is the origin and that the tangent hyperplane to $\partial\Omega$ at $a$ is given by $x_n = 0$. In a small neighbourhood of the origin, $\Omega$ and $\partial\Omega$ will then be given by

(54) $$\Omega = \{x : x_n > \phi(x')\}, \qquad \partial\Omega = \{x : x_n = \phi(x')\},$$

for some smooth function $\phi$ on $\mathbf{R}^{n-1}$ satisfying $\phi(0) = 0$, $\nabla\phi(0) = 0$. In the same neighbourhood, we can therefore use the simple map

$$\Psi(x) = (x', x_n - \phi(x'))$$

as a local chart mapping (a piece of) $\overline{\Omega}$ onto (a piece of) $\overline{\mathbf{R}^n_+}$. The Laplace operator $\Delta$ on $\Omega$ transforms under $\Psi$ into the second order operator

$$\mathcal{L} = (1 + |\nabla\phi|^2)\partial_{nn} - (\Delta\phi)\partial_n - 2\sum_{k\neq n}(\partial_k\phi)\partial_{kn} + \sum_{k\neq n}\partial_{kk}$$
$$= (1 + |\nabla\phi|^2)\partial_{nn} - (\Delta\phi)\partial_n - 2\nabla\phi \cdot \nabla'\partial_n + \Delta'$$

on $\mathbf{R}^n_+$, where we write for brevity $\partial_k = \frac{\partial}{\partial x_k}$ and $\partial_{jk} = \partial_j\partial_k$, and $\Delta', \nabla'$ denote the Laplacian and the gradient with respect to $x'$ (the prime being omitted if there is no danger of confusion).

Let us now look at the potential operator $K$ which solves the transformed problem,

$$\mathcal{L}f = 0 \text{ on } \overline{\mathbf{R}^n_+}, \quad f(\cdot, 0) = u \text{ on } \mathbf{R}^{n-1}.$$

We know from Section 3 that $K$ is of the form (17) with $k \in S^{d-1}_{\text{cl}}(\mathbf{R}^{n-1}, \mathbf{R}^{n-1}, \mathcal{S}_+)$ where $d = 0$. The boundary condition clearly forces (cf. (24)) $k(x', 0, \xi') = 1$ for all $x'$ and $\xi'$, or, in terms of the homogeneous components (50),

(55) $$k_0(x', 0, \xi') = 1, \quad k_{-j}(x', 0, \xi') = 0 \text{ for } j > 0, \qquad \forall x', \xi' \in \mathbf{R}^{n-1}.$$

On the other hand, $\mathcal{L}K = 0$ means that

(56) $$\mathcal{L}[e^{ix'\cdot\xi'}k(x', x_n, \xi')] = 0 \qquad \forall x', x_n, \xi'$$

($\mathcal{L}$ is being applied to the $x', x_n$ variables). Now

$$e^{-ix'\cdot\xi'}\mathcal{L}[e^{ix'\cdot\xi'}k(x', x_n, \xi')] = (1+|\nabla\phi|^2)\partial_{nn}k - (\Delta\phi)\partial_n k - 2\nabla\phi \cdot (i\xi'\partial_n k + \nabla'\partial_n k)$$
$$+ (\Delta'k + 2i\xi' \cdot \nabla'k - |\xi'|^2 k).$$

Substituting for $k$ the polyhomogeneous expansion $\sum_j k_{-j}$ and comparing terms with the same degrees of homogeneity, (56) yields the system of equations

(57) $$\mathcal{M}k_{-j} = \mathcal{R}k_{-j+1} + \mathcal{Z}k_{-j+2}, \quad j = 0, 1, 2, \ldots, \quad k_1 = k_2 := 0,$$

where $\mathcal{M}, \mathcal{R}, \mathcal{Z}$ are ordinary differential operators in $x_n$:

$$\mathcal{M}g = -(1 + |\nabla\phi|^2)\partial_{nn}g + 2i(\xi' \cdot \nabla\phi)\partial_n g + |\xi'|^2 g,$$
$$\mathcal{R}g = -(\Delta\phi)\partial_n g - 2\nabla\phi \cdot \nabla'\partial_n g + 2i(\xi' \cdot \nabla' g),$$
$$\mathcal{Z}g = \Delta' g.$$

Also, $k_{-j}(x', x_n, \xi')$ need to decay rapidly as $x_n \to +\infty$ for each fixed $x', \xi'$ (because they belong to $\mathcal{S}_+$ with respect to $x_n$). From this it is possible to solve (57), with the boundary conditions (55), recursively with a unique solution at each step and obtain all the $k_{-j}$.

For $j = 0$, (57) just reads $\mathcal{M}k_0 = 0$, so $k_0$ is a linear combination of $e^{ix_n \eta_+(x',\xi')}$ and $e^{ix_n \eta_-(x',\xi')}$, where

$$\eta_\pm(x', \xi') = \frac{\xi' \cdot \nabla\phi \pm i\sqrt{|\xi'|^2(1 + |\nabla\phi|^2) - |\xi' \cdot \nabla\phi|^2}}{1 + |\nabla\phi|^2}.$$

The requirement of rapid decay in $x_n$ means that only $\eta_+$ enters, and (55) then yields

$$k_0(x', x_n, \xi') = e^{ix_n \eta_+(x',\xi')}.$$

Next, for $j = 1$, (57) reads

(58) $$\mathcal{M}k_{-1} = \mathcal{R}k_0.$$

Note that $\eta_+$ is homogeneous in $\xi'$ of degree 1. A moment's computation therefore reveals that the right-hand side of (58) has the form

$$\mathcal{R}k_0 = F_1(x', \xi')k_0 + x_n F_2(x', \xi')k_0$$

where $F_1, F_2$ are homogeneous in $\xi'$ of degree 1 and 2, respectively. Solving (58), one finds that $k_{-1}$ has to be of the form $G(x', \xi')k_0 + G_0(x', \xi')x_n k_0 + G_1(x', \xi')x_n^2 k_0$. The boundary condition (55) forces $G \equiv 0$, and thus

$$k_{-1}(x', x_n, \xi') = [G_0(x', \xi')x_n + G_1(x', \xi')x_n^2]e^{ix_n \eta_+(x',\xi')},$$

with $G_0, G_1$ homogeneous in $\xi'$ of degree 0 and 1, respectively. Continuing in this way, it transpires that for all $j > 0$,

(59) $$k_{-j}(x', x_n, \xi') = \sum_{q=1}^{2j} G_{qj}(x', \xi')x_n^q e^{ix_n \eta_+(x',\xi')},$$

for some $G_{qj}(x', \xi')$ homogeneous in $\xi'$ of degree $q - j$.

The contribution from $k_{-j}$ to the Schwartz kernel $k_K(x, y)$ of $K$ (i.e. to our Poisson kernel $K(x, \zeta)$ in the local chart $\Psi$) is thus given by the inverse Fourier transform with respect to $\xi'$ of

$$\sum_{q=1}^{2j} x_n^q G_{qj}(x', \xi')\theta(\xi')e^{ix_n \eta_+(x',\xi')}$$

evaluated at $x' - y'$. (For $j = 0$, this is to be interpreted just as $\theta(\xi')e^{ix_n\eta_+(x',\xi')}$.) In particular, at the point $x' = 0$ that we are interested in, we have $\nabla\phi(0) = 0$ so $\eta_+(0,\xi') = i|\xi'|$ and we need the inverse Fourier transform with respect to $\xi'$ of

$$(60) \qquad \sum_{q=1}^{2j} x_n^q G_{qj}(\xi')\theta(\xi')e^{-x_n|\xi'|} =: \sum_{q=1}^{2j} x_n^q \theta(\xi') F_{qj,x_n}(\xi')$$

(we omit the argument $x' = 0$ in $G_{qj}$ for brevity). Analogously as in §4.2, one shows that for $q - j > 1 - n$, i.e. when $F_{qj,x_n}(\xi') = G_{qj}(\xi')e^{-x_n|\xi'|}$ is integrable at the origin, $\check{F}_{qj,x_n}(x') = x_n^{1-n-q+j}\check{F}_{qj,1}(\frac{x'}{x_n})$ is homogeneous on $\mathbf{R}_+^n$ of degree $1-n-q+j$; secondly, as $F_{qj,1} \in \mathcal{C}_{q-j}$, whence $\check{F}_{qj,1} \in S^{1-n-q+j}$, that $\check{F}_{qj,x_n}(x')$ is smooth in $(x', x_n) \in \mathbf{R}_+^n$ and extends smoothly to $\overline{\mathbf{R}_+^n} \setminus \{0\}$; thirdly, that inserting the patch function $\theta$ introduces only an error term smooth on all of $\overline{\mathbf{R}_+^n}$; and finally, differentiating and then integrating back $m$ times with respect to $x_n$, where $m = [2-n-q+j]$, that all the facts just mentioned remain in force also for $q-j \leq 1-n$. Thus the inverse Fourier transform of (60) with respect to $\xi'$ has the form

$$(61) \qquad \sum_{q=1}^{2j} x_n^q \big[U_{1-n-q+j}(x) + Q_{1-n-q+j}(x)\log|x|\big] + s(x)$$

with $U_{1-n-q+j}, Q_{1-n-q+j}$ smooth and homogeneous of the indicated degrees on $\overline{\mathbf{R}_+^n} \setminus \{0\}$, with $Q_{1-n-q+j}$ a polynomial, and $s \in C^\infty(\mathbf{R}^n)$. Again, for $j = 0$, the sum in (61) has to be replaced just by the inverse Fourier transform with respect to $\xi'$ of $e^{-x_n|\xi'|}$, which is well known to be equal to (see e.g. [22], p. 247)

$$(62) \qquad \frac{c_n x_n}{(|x'|^2 + x_n^2)^{n/2}}, \qquad c_n = \frac{\Gamma(\frac{n}{2})}{\pi^{n/2}}.$$

For $j > 0$, the log-term in (61) appears only for $1 - n - q + j \geq 0$; since $q \geq 1$, it therefore occurs only for $j \geq n$. Summing over $j$, it thus transpires that the Poisson kernel $K(x,\zeta) = \check{k}(x', x_n, x' - y')$ in local chart satisfies

$$\check{k}(y', x_n, x') = \frac{c_n x_n}{|x|^n} F\Big(y', |x|, \frac{x}{|x|}\Big) + x_n G(y', x)\log|x| + s(y', x),$$

with $c_n$ as in (62), $F \in C^\infty(\mathbf{R}^{n-1} \times \overline{\mathbf{R}_+} \times \mathbf{S}^{n-1})$ and $G, s \in C^\infty(\mathbf{R}^{n-1} \times \overline{\mathbf{R}_+^n})$. Since, from its very definition, $K(x,y')$ vanishes for $x_n = 0$ (and $x' \neq y'$), it must be possible to pull out a factor of $x_n$ also from $s(y', x)$, and then merge the latter into the first summand. Altogether, we have thus arrived at the following more precise version of Theorem 8 for the particular case of the Poisson kernel, which we have stated as Theorem 1 in the Introduction.

**Theorem 13.** (= Theorem 1) *Let $\Omega$ be a bounded domain in $\mathbf{R}^n$ with smooth boundary. Then the Poisson kernel $K(x,\zeta)$ for the Dirichlet problem (2) can be written in the form*

$$K(x,\zeta) = \frac{c_n d(x)}{|x-\zeta|^n}\bigg[F\Big(\zeta, |x-\zeta|, \frac{x-\zeta}{|x-\zeta|}\Big) + G(x,\zeta)|x-\zeta|^n \log|x-\zeta|\bigg]$$

where $d \in C^\infty(\overline{\Omega})$, $d > 0$ on $\Omega$ and $d(x) = \text{dist}(x, \partial\Omega)$ for $x$ near $\partial\Omega$, $c_n$ is as in (62), $F \in C^\infty(\partial\Omega \times \overline{\mathbf{R}_+} \times \mathbf{S}^{n-1})$, $F(\cdot, 0, \cdot) = 1$, and $G \in C^\infty(\overline{\Omega} \times \partial\Omega)$.

We conclude this section by showing that the log term is indeed present in general, although the explicit formulas

$$K(x, \zeta) = \frac{c_n}{2} \frac{1 - |x|^2}{|x - \zeta|^n}$$

for $\Omega = \{x \in \mathbf{R}^n : |x| < 1\}$, and

$$K(x, \zeta) = \frac{c_n x_n}{|x - \zeta|^n}$$

for $\Omega = \mathbf{R}^n_+$, show that it is absent for the ball and the half-space. To that end, consider the special domain

$$\Omega = \{y \in \mathbf{R}^n : y_n > \Phi(|y'|^2)\}$$

for some smooth function $\Phi$ on $\overline{\mathbf{R}_+}$; i.e. the special case of (54) where $\phi(y') = \Phi(|y'|^2)$ depends only on $|y'|$. The computations above then simplify somewhat; for instance, one has the explicit inverse to the operator $\mathcal{M}$, namely, the unique solution to

$$\mathcal{M}g = x_n^m e^{ix_n \eta_+}$$

($m = 0, 1, 2, \dots$), rapidly decreasing as $x_n \to +\infty$, is

$$g = \frac{m!}{(-\upsilon)^{m+2}} \sum_{j=1}^{m+1} \frac{(-x_n \upsilon)^j}{j! i^{m+2-j}} \frac{e^{ix_n \eta_+}}{1 + |\nabla \phi|^2},$$

where $\upsilon := \eta_+ - \eta_-$. Using all this, one arrives at the following formulas (checked on a computer) for $n = 3$ and $x' = 0$:

$$k_0(0, x_n, \xi') = e^{-x_n |\xi'|},$$

$$k_{-1}(0, x_n, \xi') = x_n(1 - x_n |\xi'|) \Phi'(0) e^{-x_n |\xi'|},$$

$$k_{-2}(0, x_n, \xi') = \left(\frac{x_n}{2|\xi'|} + \frac{5}{2} x_n^2 - 3|\xi'| x_n^3 + \frac{|\xi'|^2 x_n^4}{2}\right) \Phi'(0)^2 e^{-x_n |\xi'|},$$

$$k_{-3}(0, x_n, \xi') = \left(x_n \frac{2\Phi'(0)^3 - \Phi''(0)}{2|\xi'|^2} + x_n^2 \frac{4\Phi'(0)^3 - \Phi''(0)}{2|\xi'|} + x_n^3 (7\Phi'(0)^3 - \Phi''(0))\right.$$
$$\left. - \frac{19\Phi'(0)^3 - \Phi''(0)}{2} x_n^4 |\xi'| + \frac{5\Phi'(0)^2}{2} x_n^5 |\xi'|^2 - \frac{\Phi'(0)^3}{6} x_n^6 |\xi'|^3\right) e^{-x_n |\xi'|},$$

and, hence, after a small computation[2]

(63)
$$K((0, x_n), x') = \frac{t}{2\pi |x|^2} + \frac{t^2(2 - 3t^2)\Phi'(0)}{2\pi |x|}$$
$$+ \frac{t(1 - t^2)(1 + 12t^2 - 15t^4)\Phi'(0)^2}{4\pi}$$
$$+ \frac{t(\Phi''(0) - 2\Phi'(0)^3)}{4\pi} |x| \log |x| + \text{(higher order terms)},$$

---

[2] See (70) below for a hint.

where we have set for brevity $t = x_n/|x|$. (All these formulas relate to the local chart by $\overline{\mathbf{R}^3_+}$ obtained from $\Omega$ via the map $(y', y_n) \mapsto (y', y_n - \Phi(|y'|^2)) \equiv (x', x_n)$.) Thus the log term is present whenever $\Phi''(0) - 2\Phi'(0)^3 \neq 0$. Note that the last expression vanishes if $\Omega$ is the half-space (i.e. $\Phi \equiv 0$), or if $\Omega$ is the ball of radius $R$ tangent to $y_n = 0$ at the origin (i.e. $\Phi(t) = R - \sqrt{R^2 - t}$, so that $\Phi'(0) = \frac{1}{2R}$, $\Phi''(0) = \frac{1}{4R^3}$), in complete accordance with the fact that there is no log term in $K(x, \zeta)$ in these cases.

Finally, we remark that the log-term is always absent in dimension $n = 2$. Indeed, since harmonic functions in $\mathbf{C} \cong \mathbf{R}^2$ are preserved by composition with holomorphic maps (this is no longer true in $\mathbf{C}^m \cong \mathbf{R}^{2m}$ when $m > 1$, nor has any good analogue in odd dimensions), while any simply connected domain in $\mathbf{C}$ with smooth boundary can be biholomorphically mapped onto the disc (by the Riemann mapping theorem), and the Poisson kernel of the disc has no log-term, the Poisson kernel of the original domain has no log-term either. (The argument extends in fact also to multiply connected domains with smooth boundary, in view of the local character of the boundary $\Psi$DOs.) Thus the examples above for $n = 3$ are probably the simplest domains for which the log-term is present.

## 6. The harmonic Bergman kernel

The singularity of the Schwartz kernel of a singular Green operator $G$ can again be obtained from the asymptotic expansion of its symbol

$$(64) \qquad g \sim \sum_{j=0}^{\infty} g_{d-j}$$

into the homogeneous components satisfying

$$(65) \qquad g_{d-j}(x', \tfrac{x_n}{\lambda}, \tfrac{y_n}{\lambda}, \lambda \xi') = \lambda^{d-j+1} g_{d-j}(x', x_n, y_n, \xi')$$

for $\lambda \geq 1$, $|\xi'| \geq 1$; as with potential operators, we will again assume from now on that (65) actually holds for all $\lambda > 0$ and $\xi' \neq 0$, at the expense of having $g_{d-j}$ singular at $\xi' = 0$ and the asymptotic expansion (64) holding only for $|\xi'| \geq 1$, analogously as in (52). The contribution to the Schwartz kernel

$$(66) \qquad k_G(x, y) = \check{g}(x', x_n, y_n, z'), \qquad z' = x' - y',$$

of $G$ from $g_{d-j}$ is again given by the inverse Fourier transform with respect to $\xi'$ of

$$g_{d-j}(x', x_n, y_n, \xi')\theta(\xi')$$

for some patch function $\theta$. To obtain the $g_{d-j}$ more explicitly, we use the result for the Poisson operator from Section 5 together with the composition rules reviewed in §2.3.

Recall that our singular Green operator of interest, the harmonic Bergman projection, is given by (31)

$$G = \mathbf{K}\Lambda^{-1}\mathbf{K}^*, \qquad \Lambda := \mathbf{K}^*\mathbf{K},$$

where $\mathbf{K}$ is the Poisson extension operator from Section 5. By the general rules for the Boutet de Monvel calculus, $\Lambda$ is a $\Psi$DO on the boundary of order $-1$, and $G$ is a singular Green operator of order 0 and class zero. Let us now compute the symbol of $G$. In terms of the symbol $k(x', x_n, \xi')$ of $\mathbf{K}$ (in a local chart by $\overline{\mathbf{R}^n_+}$), the adjoint $\mathbf{K}^*$ is a trace operator of class zero with symbol $k^*(x', x_n, \xi')$ given by (cf. (25))

$$k^*(x', x_n, \xi') \sim \sum_\alpha \frac{1}{\alpha!} \partial^\alpha_{x'} D^\alpha_{\xi'} \overline{k(x', x_n, \xi')}.$$

This formula is well adapted with respect to the homogeneous expansion (64): namely, one gets (remember that $\mathbf{K}$ is of degree $d = 0$) $k^* \sim \sum_{j=0}^\infty k^*_{-j}$, with

$$k^*_{-j} = \sum_{m+|\alpha|=j} \frac{1}{\alpha!} \partial^\alpha_{x'} D^\alpha_{\xi'} \overline{k_{-m}}.$$

We have seen that $k_{-j}$ are of the form (59) (where for $j = 0$ the sum $\sum_{q=1}^{2j} G_{qj} x_n^q$ has to be interpreted as 1). The last formula and a quick check show that $k^*_{-j}$ are of the same form, only with $e^{ix_n \eta_+}$ replaced by $\overline{e^{ix_n \eta_+}} = e^{-ix_n \eta_-}$; in particular, $k^*_0 = e^{-ix_n \eta_-}$.

Next, $\Lambda = \mathbf{K}^* \mathbf{K}$ is a $\Psi$DO on $\partial\Omega$ ($= \mathbf{R}^{n-1}$ in the local chart) with symbol given by (22)

$$s(x', \xi') = \sum_\alpha \frac{1}{\alpha!} \int_0^\infty D^\alpha_{\xi'} k^*(x', x_n, \xi') \partial^\alpha_{x'} k(x', x_n, \xi') \, dx_n.$$

This is again well-behaved with respect to the polyhomogeneous grading, and using the fact that

$$\int_0^\infty x_n^k e^{ix_n \eta_+ - ix_n \eta_-} \, dx_n = \frac{k!}{(i\eta_- - i\eta_+)^{k+1}},$$

we get $s \sim \sum_{j=0}^\infty s_{-j-1}$, with

$$s_{-m-1} = \sum_{l+j+|\alpha|=m} \frac{1}{\alpha!} \int_0^\infty D^\alpha_{\xi'} k^*_{-j} \partial^\alpha_{x'} k_{-l} \, dx_n,$$

and, in particular, $s_{-1} = 1/(i\eta_- - i\eta_+)$. The symbol $p$ of the inverse $\Lambda^{-1}$ therefore has the asymptotic expansion $p \sim \sum_{j=0}^\infty p_{1-j}$, given recursively by

$$p_1 = \frac{1}{s_{-1}}, \quad p_{1-m} = -\frac{1}{s_{-1}} \sum_{k+j+|\alpha|=m, j<m} \frac{1}{\alpha!} D^\alpha_{\xi'} p_{1-j} \partial^\alpha_{x'} s_{-1-k} \quad \text{for } m > 0,$$

by (15); in particular, $p_1 = i(\eta_- - \eta_+)$.

By (22) again, we get that $\mathbf{K}\Lambda^{-1}$ is a potential operator of degree 1 with symbol $v \sim \sum_{j=0}^\infty v_{1-j}$ with

$$v_{1-m} = \sum_{j+|\alpha|+l=m} \frac{1}{\alpha!} D^\alpha_{\xi'} k_{-j} \partial^\alpha_{x'} p_{1-l}.$$

For the terms with $l < m$, (59) and a routine check imply that they are again of the form (59) (with $m$ in the place of $j$) multiplied by a function homogeneous in $\xi'$ of degree $1 - l$; while for $l = m$, the corresponding term is just $p_{1-m}e^{ix_n\eta_+}$. Altogether, we thus see that $v_{1-m}$ has the form

$$v_{1-m} = \sum_{q=0}^{2m} F_{mq}(x', \xi') x_n^q e^{ix_n\eta_+} \tag{67}$$

with $F_{mq}$ homogeneous in $\xi'$ of degree $1 - m + q$ and $F_{m0} = p_{1-m}$. (Note that (67) is similar to (59), except that the sum starts from $q = 0$.) In particular, $v_1 = i(\eta_- - \eta_+)e^{ix_n\eta_+}$.

Finally, using (22) one more time, we get that $\mathbf{K}\Lambda^{-1}\mathbf{K}^* = G$ has symbol $g \sim \sum_{j=0}^{\infty} g_{1-j}$, with

$$g_{1-m} = \sum_{j+l+|\alpha|=m} \frac{1}{\alpha!} D_{\xi'}^{\alpha} v_{1-j}(x', x_n, \xi') \partial_{x'}^{\alpha} k_{-l}^*(x', y_n, \xi').$$

In particular, $g_1 = i(\eta_- - \eta_+)e^{ix_n\eta_+ - iy_n\eta_-}$, while for general $m$ we again have from (59) (for $k^*$) and (67)

$$g_{1-m} = \sum_{r,q=0}^{2m} F_{mrq}(x', \xi') x_n^q y_n^r e^{i(x_n\eta_+ - y_n\eta_-)}, \tag{68}$$

with $F_{mrq}$ homogeneous of degree $1 - m + r + q$ in $\xi'$ and $F_{m00} = p_{1-m}$.

The contribution to the singularity of the Schwartz kernel (66) of $G$ — i.e. to the harmonic Bergman kernel $H$ — is given by the inverse Fourier transform with respect to $\xi'$, evaluated at $x' - y'$, of (68) multiplied by a patch function $\theta(\xi')$. Going through the proof in §4.3 (analogously to the argument after (60) in the preceding section), it transpires that the latter equals, modulo smooth functions on $\mathbf{R}^{n-1} \times \mathbf{R}^{n+1}$,

$$\begin{aligned}x_n^q y_n^r \Big[ U_{-n+m-r-q}(x', z) + Q_{-n+m-r-q}(x', z) \log |z| \Big], \\ z = (x' - y', x_n, y_n) \in \mathbf{R}^{n+1},\end{aligned} \tag{69}$$

for some $U_{-n+m-r-q}(x', z), Q_{-n+m-r-q}(x', z)$ smooth on $\mathbf{R}^{n-1} \times (\overline{\mathbf{R}_+^{n+1}} \setminus \{0\})$ and homogeneous in $z$ of the indicated degree, with $Q$ a polynomial in $z$.

In particular, at $x' = 0$, we have $i\eta_- = -i\eta_+ = |\xi'|$, so $g_1 = 2|\xi'|e^{-(x_n+y_n)|\xi'|}$, and

$$\check{g}_1 = -\Big[2\frac{\partial}{\partial t}e^{-t|\xi'|}\Big]^{\vee}\Big|_{t=x_n+y_n} = -2\frac{\partial}{\partial t}\frac{c_n t}{(t^2 + |x'|^2)^{n/2}}\Big|_{t=x_n+y_n}, \tag{70}$$

with $c_n$ as in (62), which shows that the leading term of $k_G(x, y)$ in the local chart is

$$2c_n \frac{(n-1)(x_n + y_n)^2 - |x' - y'|^2}{[(x_n + y_n)^2 + |x' - y'|^2]^{n/2+1}}.$$

In terms of the boundary distance $d(x)$ and the reflected point $\tilde{y} = (y', -y_n)$, this becomes

$$2c_n \frac{n(d(x) + d(y))^2 - |x - \tilde{y}|^2}{|x - \tilde{y}|^{n+2}}.$$

Since the differential of our local chart $\Psi$ from Section 5 equals the identity at the origin, the last expression remains in force also back in $\Omega$. We have thus arrived at the following slight strengthening of Corollary 12 for the particular case of the singular Green operator $G = \mathbf{K}\Lambda^{-1}\mathbf{K}^*$ of order 0 and class zero, which has been stated as Theorem 2 in the Introduction.

**Theorem 14.** (= Theorem 2) *Let $\Omega$ be a bounded domain in $\mathbf{R}^n$ with smooth boundary. The harmonic Bergman kernel $H(x,y)$ of $\Omega$ is $C^\infty$ for $(x,y)$ away from the boundary diagonal, while near $\operatorname{diag}\partial\Omega$ it can be written in the form*

$$(71) \qquad H(x,y) = \frac{2c_n}{|x-\tilde{y}|^n} F\Big(x, y, |x-\tilde{y}|, \frac{x-\tilde{y}}{|x-\tilde{y}|}\Big) + G(x,y) \log|x-\tilde{y}|,$$

*where $G \in C^\infty(\overline{\Omega} \times \overline{\Omega})$, $F \in C^\infty(\overline{\Omega} \times \overline{\Omega} \times \overline{\mathbf{R}_+} \times \mathbf{S}^{n-1})$, $F(x,x,0,\nu) = n\langle \nu, \nabla d(x)\rangle^2 - 1$ for $x \in \partial\Omega$, and $c_n$ is as in (62).*

Note that, in particular, the leading term of the restriction $H|_{\partial\Omega\times\partial\Omega}$ of $H$ to the boundary is simply $-2c_n|x-y|^{-n}$.

*Remark.* In order to obtain an expression for $H(x,y)$ like (71), valid on all of $\Omega\times\Omega$, not only near the boundary, one possibility is as follows:

$$(72) \qquad H(x,y) = \frac{2c_n}{|v(x,y)|^n} \Big[ F\Big(x, y, |v(x,y)|, \frac{v(x,y)}{|v(x,y)|}\Big) + G(x,y) \log|v(x,y)| \Big],$$

where $G \in C^\infty(\overline{\Omega} \times \overline{\Omega})$, $F \in C^\infty(\overline{\Omega} \times \overline{\Omega} \times \overline{\mathbf{R}_+} \times \mathbf{S}^{n+1})$, and $v : \overline{\Omega} \times \overline{\Omega} \to \mathbf{R}^{n+1}$ is any smooth function such that

$$v(x,y) = (\rho(x) - \rho(y), d(x), d(y))$$

for $(x,y)$ in some small neighbourhood of $\operatorname{diag}\partial\Omega$, where $\rho(x)$ is the point of $\partial\Omega$ closest to $x$ ($\rho$ is well defined in a small neighbourhood of $\partial\Omega$, and such $v$ exists by Seeley's extension theorem). One can further replace $|v(x,y)|$ in (72) by $d_\Delta(x,y)$, with any $d_\Delta \in C^\infty(\overline{\Omega} \times \overline{\Omega})$ such that $d_\Delta > 0$ on $\overline{\Omega \times \Omega} \setminus \operatorname{diag}\partial\Omega$ and $d_\Delta(x,y) = \operatorname{dist}((x,y), \operatorname{diag}\partial\Omega)$ for $(x,y)$ near $\operatorname{diag}\partial\Omega$, etc.; cf. the remark after the proof of Theorem 10. $\square$

We conclude this section again by showing that the log term is present, i.e. $G \not\equiv 0$, in general, although this is not the case for the familiar situations like the ball $\Omega = \{x \in \mathbf{R}^n : |x| < 1\}$ with

$$H(x,y) = \frac{c_n}{2} \frac{(n-4)|x|^4|y|^4 + (8x\cdot y - 2n - 4)|x|^2|y|^2 + n}{(1 - 2x\cdot y + |x|^2|y|^2)^{n/2+1}},$$

or the upper half-space $\mathbf{R}^n_+$ with

$$H(x,y) = 2c_n \frac{n(x_n + y_n)^2 - |x-\tilde{y}|^2}{|x-\tilde{y}|^{n+2}}.$$

To this end, we use once more the domains from Section 5

$$\Omega = \{y \in \mathbf{R}^n : y_n > \Phi(|y'|^2)\}$$

with some $\Phi \in C^\infty(\overline{\mathbf{R}_+})$ and $n = 3$. We again pass to the local chart $\overline{\mathbf{R}^n_+}$ by means of the map $\Psi : (y', y_n) \mapsto (y', y_n - \Phi(|y'|^2)) = (x', x_n)$. We will consider $H(x,y)$ only at $x' = 0$ and $x_n = y_n = 0$ (that is, we exhibit a log-term even among the terms (69) with $q = r = 0$). Now since $F_{m00} = p_{1-m}$, we have that, modulo smooth

POISSON AND HARMONIC BERGMAN KERNELS 35

functions, $H(0,(y',0)) = \check{g}(0,0,0,-y')$ is the inverse Fourier transform with respect to $\xi'$, evaluated at $-y'$, of

$$g(0,0,0,\xi') = p(0,\xi').$$

Combining the formulas for $k^*$, $s$ and $p$ above with those for $k_0$, $k_{-1}$, $k_{-2}$, $k_{-3}$ in Section 5, we obtain after some calculation (again verified on a computer)

$$s(0,\xi') \sim \frac{1}{2|\xi'|} - \frac{\Phi'(0)}{|\xi'|^2} + \frac{5\Phi'(0)^2}{4|\xi'|^3} + \frac{5\Phi''(0)}{2|\xi'|^4} + \dots$$

and

$$p(0,\xi') \sim 2|\xi'| + 4\Phi'(0) - \frac{\Phi'(0)^2}{|\xi|} + \frac{2\Phi''(0) - 4\Phi'(0)^3}{|\xi'|^2} + \dots.$$

Consequently,

$$(73) \quad \check{p}(0,y') = -\frac{2c_3}{|y'|^3} - \frac{c_3\Phi'(0)^2}{|y'|} - c_3(2\Phi''(0) - 4\Phi'(0)^3)\log|y'|$$
$$+ \text{(higher order terms)}.$$

Hence the log-term is nonzero as soon as $2\Phi''(0) - 4\Phi'(0)^3 \neq 0$. Note that, again, the last expression vanishes for the upper half-space ($\Phi \equiv 0$) as well as for the ball ($\Phi(t) = R - \sqrt{R^2 - t}$, $\Phi'(0) = \frac{1}{2R}$, $\Phi''(0) = \frac{1}{4R^3}$), as it should.

## 7. Concluding remarks

**7.1 Weighted harmonic Bergman kernels.** For a positive weight function $w$ on $\Omega$, one can consider the subspace $L^2_{\text{harm}}(\Omega, w)$ of harmonic functions in the weighted space $L^2(\Omega, w)$; if $w^{-1}$ is locally bounded (in particular, if $w$ is continuous), the standard argument shows that point evaluations are continuous on $L^2_{\text{harm}}(\Omega, w)$, and uniformly so on compact subsets, so that $L^2_{\text{harm}}(\Omega, w)$ is actually closed in $L^2(\Omega, w)$ and possesses a reproducing kernel — the weighted harmonic Bergman kernel $H_w(x,y)$. Of course, for $H_w$ not to be identically zero, the space $L^2_{\text{harm}}(\Omega, w)$ needs to be nontrivial; this is the case e.g. whenever $w$ is integrable over $\Omega$ (since the constant functions then belong to $L^2_{\text{harm}}(\Omega, w)$).

Much of what has been said above for the unweighted case remains also true in this weighted situation; in particular, the adjoint of the Poisson operator $\mathbf{K}$ with respect to the weighted inner product on $\Omega$ is just $\mathbf{K}^* w$ (i.e. the unweighted adjoint $\mathbf{K}^*$ preceded by multiplication by $w$), and the operator

$$(74) \quad \Lambda_w := \mathbf{K}^* w \mathbf{K}$$

is related to the weighted Bergman kernel by

$$(75) \quad H_w(x,y) = \langle \Lambda_w^{-1} K_x, K_y \rangle_{\partial\Omega}, \qquad K_x \equiv K(x,\cdot),$$

and to the weighted harmonic Bergman projection $\mathbf{\Pi}_w$ by

$$\mathbf{\Pi}_w = \mathbf{K} \Lambda_w^{-1} \mathbf{K}^* w.$$

The Schwartz kernel of $\mathbf{\Pi}_w$, however, is not equal to $H_w$, since the Schwartz kernel is taken "with respect to the Lebesgue measure", while $H_w$ must be integrated against the weight $w$; for this reason, $H_w$ is actually the Schwartz kernel of

$$G := \mathbf{\Pi}_w w^{-1} = \mathbf{K} \Lambda_w^{-1} \mathbf{K}^*.$$

In particular, if $\Omega$ is bounded with smooth boundary and $w$ is of the form

(76) $$w(x) = d(x)^\alpha e^{g(x)},$$

with $d \in C^\infty(\overline{\Omega})$, $d > 0$ on $\Omega$ and $d(x) = \text{dist}(x, \partial\Omega)$ for $x$ near $\partial\Omega$, $\alpha > -1$ (to make $w$ integrable), and $g \in C^\infty(\overline{\Omega})$, then the operator (74) belongs to the Boutet de Monvel calculus: namely, $\Lambda_w$ is a $\Psi$DO on $\partial\Omega$ of order $-\alpha - 1$ and with nonvanishing principal symbol (see [8])

(77) $$\frac{\Gamma(\alpha+1) e^{g(x')}}{(2|\xi'|)^{\alpha+1}}.$$

Since $\mathbf{K}$ is a potential operator of order 0 while $\mathbf{K}^*$ is a trace operator of class zero and order $-1$, it follows that $G$ is a singular Green operator of class zero and order $\alpha$, with principal symbol (in the local chart from Section 5, at $x' = 0$)

(78) $$g_{\alpha+1}(x', x_n, y_n, \xi') = \frac{(2|\xi'|)^{\alpha+1} e^{-(x_n+y_n)|\xi'|}}{\Gamma(\alpha+1) e^{g(x')}}.$$

It follows that $H_w(x, y)$ is given by Theorem 10 (and Corollary 12) with $d = \alpha$, and the leading term (in the local chart by $\overline{\mathbf{R}^n_+}$) can be obtained as the inverse Fourier transform with respect to $\xi'$, evaluated at $x' - y'$, of (78). (For $\alpha$ integer, the latter can be computed explicitly as

$$\check{g}_{\alpha+1}(0, x_n, y_n, x') = \left(-2\frac{\partial}{\partial t}\right)^{\alpha+1} \frac{[e^{-t|\xi'|}]^\vee(x')}{\Gamma(\alpha+1) e^{g(x')}}\bigg|_{t=x_n+y_n}$$

$$= \frac{c_n}{\Gamma(\alpha+1) e^{g(x')}} \left(-2\frac{\partial}{\partial t}\right)^{\alpha+1} \frac{t}{(t^2+|x'|^2)^{n/2}}\bigg|_{t=x_n+y_n}.$$

For general $\alpha > -1$, $\check{g}_{\alpha+1}(0, x_n, y_n, x')$ is given by the hypergeometric function

$$\frac{2^\alpha \Gamma(\frac{n+\alpha}{2}) \Gamma(\frac{n+\alpha-1}{2})}{\pi^{n/2} \Gamma(\frac{n-1}{2}) \Gamma(\alpha+1) e^{g(x')} (x_n+y_n)^{n+\alpha-1}} \,_2F_1\left(\begin{array}{c}\frac{n+\alpha}{2}, \frac{n+\alpha-1}{2}\\ \frac{n-1}{2}\end{array}\bigg| -\frac{|x'|^2}{(x_n+y_n)^2}\right),$$

as can be seen by expanding $e^{ix'\cdot\xi'}$ into power series and integrating term by term.)

**7.2 Log-terms in dimension 2.** We have seen in Section 5 that the Poisson kernel of a smoothly bounded domain $\Omega \subset \mathbf{R}^2 \cong \mathbf{C}$ has no log-term, due to the fact that this is the case for the disc and the Riemann mapping theorem. Computations indicate that also the harmonic Bergman kernel $H(x, y)$ has no log-term for such domains; or, using again the Riemann mapping theorem, the weighted harmonic Bergman kernel $H_w$ on the disc with weight $w$ of the form $w = |f'|^2$, where $f$ is a conformal map, has no log term. It would be interesting to know if this is indeed true, and why. (For general weights $w$ on the disc, $H_w$ has log terms; an example is $w(z) = 2-|z|^2$, when $H_w(x,y) = 2\operatorname{Re}[(1-x\overline{y})^{-2}-(1-x\overline{y})^{-1}-2\log(1-x\overline{y})+F(x\overline{y})]$, with $F$ continuous on the closed disc.)

For our sample domains $\{x \in \mathbf{R}^3 : x_3 > \Phi(|x'|^2)\}$ in Sections 5 and 6, it was also the case that the conditions for the presence of the log-terms in $K(x,\zeta)$ and in $H(x,y)$ were the same — namely, $\Phi''(0) - 2\Phi'(0)^3 \neq 0$; suggesting that perhaps the log term is present in $K(x,\zeta)$ if and only if it is present in $H(x,y)$. The present author has no idea whether this is indeed true, nor does he know whether there are other domains than the ball and the half-space for which the log-term in either $K(x,\zeta)$ or $H(x,y)$ would be absent.

**7.3 Harmonic Sobolev-Bergman kernels.** Our methods apply also to the reproducing kernels of Sobolev type spaces of harmonic functions. For instance, consider the space

$$W^1_{\text{harm}}(\Omega) = \{u \text{ harmonic on } \Omega : u, \nabla u \in L^2(\Omega)\}$$

with the standard norm

$$\|u\|_1^2 := \|u\|_\Omega^2 + \sum_{j=1}^n \|\partial_j u\|_\Omega^2.$$

Note that the function $\partial_j u$ is harmonic if $u$ is; thus, for $f \in C^\infty(\partial\Omega)$, $\partial_j \mathbf{K} u = \mathbf{K} R_j u$, where

$$R_j = \gamma \partial_j \mathbf{K} = \Lambda^{-1} \mathbf{K}^* \partial_j \mathbf{K}$$

is a $\Psi$DO on $\partial\Omega$ of order 1, by Boutet de Monvel's calculus. Then

$$\|\mathbf{K}f\|_1^2 = \langle \mathbf{K}f, \mathbf{K}f\rangle_\Omega + \sum_{j=1}^n \langle \mathbf{K}R_j f, \mathbf{K}R_j f\rangle_\Omega = \langle Tf, f\rangle_{\partial\Omega},$$

where

$$(79) \qquad T = \mathbf{K}^*\mathbf{K} + \sum_{j=1}^n R_j^* \mathbf{K}^* \mathbf{K} R_j = \Lambda + \sum_{j=1}^n R_j^* \Lambda R_j$$

is a positive self-adjoint $\Psi$DO on $\partial\Omega$ of order 1, with principal symbol $\sigma(T) = \sigma(\Lambda)\sum_j |\sigma(R_j)|^2$. One can show that

$$\sum_{j=1}^n |\sigma(R_j)|^2 = \frac{1}{2\sigma(\Lambda)^2};$$

thus $\sigma(T) = \frac{1}{2}\sigma(\Lambda)^{-1} > 0$, so $T$ is in fact elliptic.

Consider now, quite generally, the Hilbert space of harmonic functions on $\Omega$ obtained as the completion of $\{\mathbf{K}f : f \in C^\infty(\partial\Omega)\}$ with respect to the norm

$$\|\mathbf{K}f\|_T^2 := \langle Tf, f\rangle,$$

for some positive self-adjoint elliptic $\Psi$DO $T$ on $\partial\Omega$. For the corresponding reproducing kernel $H_T(x,y) \equiv H_{T,x}(y)$, we then get

$$\langle f, K_x\rangle_{\partial\Omega} = \mathbf{K}f(x) = \langle \mathbf{K}f, H_{T,x}\rangle_T = \langle Tf, \gamma H_{T,x}\rangle_{\partial\Omega}$$

by the definition of the $T$-norm. Hence $K_x = T\gamma H_{T,x}$, or

$$\gamma H_{T,x} = T^{-1} K_x \tag{80}$$

and

$$H_T(x,y) = \mathbf{K}\gamma H_{T,x}(y) = \langle \gamma H_{T,x}, K_y \rangle_{\partial\Omega} = \langle T^{-1} K_x, K_y \rangle_{\partial\Omega},$$

generalizing the formulas (6) for $T = \mathbf{K}^*\mathbf{K}$ and (75) for $T = \mathbf{K}^*w\mathbf{K}$. Finally, $H_T$ is the Schwartz kernel of the operator $G$ given by

$$\begin{aligned} G\mathbf{K}f(x) &= \int_\Omega \mathbf{K}f(y) H_T(x,y)\, dy \\ &= \langle \mathbf{K}f, H_{T,x} \rangle_\Omega = \langle \mathbf{K}f, \mathbf{K}T^{-1}K_x \rangle_\Omega \qquad \text{by (80)} \\ &= \langle T^{-1}\mathbf{K}^*\mathbf{K}f, K_x \rangle_{\partial\Omega} = \mathbf{K}T^{-1}\mathbf{K}^*\mathbf{K}f(x), \end{aligned}$$

that is,

$$G = \mathbf{K}T^{-1}\mathbf{K}^*. \tag{81}$$

For the harmonic Sobolev space $W^1_{\text{harm}}$ above, with $T$ given by (79), the Boutet de Monvel calculus shows that $G$ is a singular Green operator of class zero and order $-2$, and the corresponding harmonic Sobolev-Bergman kernel $H_T$ is thus described by Theorem 10 and Corollary 12 with $d = -2$.

For the analogous harmonic Sobolev spaces $W^k_{\text{harm}}(\Omega)$ of higher integer order $k$, the corresponding operator $T$ becomes

$$T = \sum_{m=0}^{k} \sum_{j_1,\ldots,j_m=1}^{n} R^*_{j_1} \ldots R^*_{j_m} \Lambda R_{j_m} \ldots R_{j_1},$$

a positive self-adjoint elliptic $\Psi$DO on $\partial\Omega$ of order $2k-1$, so $G$ is a singular Green operator of order $-2k$ (and class zero) and Theorem 10 and Corollary 12 again apply.

The operator $T$ covers also the situation of weighted harmonic Bergman spaces in §7.1 (then $T = \mathbf{K}^*w\mathbf{K} = \Lambda_w$), and one can even combine the two and consider weighted harmonic Sobolev spaces with respect to weights $w$ of the form (76). This is important due to the fact that $L^2_{\text{harm}}(\Omega, w)$, with $w$ as in (76), is known to actually coincide with $W^{-\alpha/2}_{\text{harm}}(\Omega)$, and similarly $W^k_{\text{harm}}(\Omega, w)$ coincides with $W^{k-\alpha/2}_{\text{harm}}(\Omega)$ (in both cases, with equivalent norms); varying $k$ and $\alpha$, one can thus cover the whole range of $W^s_{\text{harm}}(\Omega)$ with any real $s$ (without resorting to interpolation, which is usually used to define $W^s$ for non-integer $s$). The corresponding reproducing kernels are still susceptible to the treatment as above.

**7.4 Analytic continuation.** It has been shown that for weights of the form (76), the corresponding *holomorphic* Bergman kernels $B_\alpha(x,y)$ can in fact be continued analytically in $\alpha$ to a meromorphic function in the entire complex plane [18]. The same idea — in fact, a much simpler variant of it — works also for our harmonic kernels here.

Namely, fix $d, g \in C^\infty(\overline{\Omega})$, with $d > 0$ on $\Omega$ and $d(x) = \text{dist}(x, \partial\Omega)$ near $\partial\Omega$, and let

$$w_\alpha(x) := d(x)^\alpha e^{g(x)}, \qquad \alpha > -1. \tag{82}$$

The associated weighted harmonic Bergman kernels $H_{w_\alpha} \equiv H_\alpha$ are then given by the formula (75)
$$H_\alpha(x,y) = \langle \Lambda_\alpha^{-1} K_x, K_y \rangle,$$
where
$$\Lambda_\alpha := \Lambda_{w_\alpha} = \mathbf{K}^* w_\alpha \mathbf{K}$$
are positive self-adjoint elliptic $\Psi$DOs on $\partial\Omega$ of order $-\alpha - 1$. Now, first of all, Proposition 10 and Remark 11 in [18] assert that $\Lambda_\alpha/\Gamma(\alpha+1)$ can be analytically continued to all $\alpha \in \mathbf{C}$ as a holomorphic family of $\Psi$DOs (which we still denote $\Lambda_\alpha$) of order $-\alpha - 1$ (with the principal symbols still given by (77), only with the factor $\Gamma(\alpha + 1)$ omitted). It therefore remains only to deal with the invertibility of $\Lambda_\alpha$. On the other hand, $\Lambda = \mathbf{K}^*\mathbf{K}$ is a positive self-adjoint operator, hence we can form its complex power $\Lambda^\alpha$ for any $\alpha \in \mathbf{C}$ by the spectral theorem; as $\Lambda$ is an elliptic $\Psi$DO of order $-1$, the familiar theory of Seeley [39] tells us that $\Lambda^\alpha$ is actually a $\Psi$DO of order $-\alpha$ with principal symbol $\sigma(\Lambda)^\alpha = (2|\xi'|)^{-\alpha}$. It follows that
$$F(\alpha) := \Gamma(\alpha+1)^{-1} \Lambda^{-\alpha} \Lambda_0^{-1} \Lambda_\alpha, \qquad \alpha \in \mathbf{C},$$
is a family of $\Psi$DOs on $\partial\Omega$, holomorphic in the usual sense (see §2.3 in [18]), of order 0, and with principal symbol equal to (cf. (77))
$$\frac{1}{\Gamma(\alpha+1)}(2|\xi'|)^\alpha \frac{2|\xi'|}{e^{g(x')}} \frac{\Gamma(\alpha+1)e^{g(x')}}{(2|\xi'|)^{\alpha+1}} = 1.$$
Thus $F(\alpha) - I$ is a holomorphic family of $\Psi$DOs of order $-1$, hence, as $\partial\Omega$ is compact, of compact operators; furthermore, $F(0) = \Lambda_0^{-1}\Lambda_0 = I$. A theorem of Gohberg ([24], Chapter I, Theorem 5.1) implies that there exists a set $U$ of isolated points in $\mathbf{C}$ such that $F(\alpha)$ is invertible for $\alpha \notin U$, and $F(\alpha)^{-1}$ is holomorphic on $\mathbf{C} \setminus U$ with poles with finite-rank residues at points of $U$. Consequently,
$$Y(\alpha) := \Gamma(\alpha+1)^{-1} F(\alpha)^{-1} \Lambda^{-\alpha} \Lambda_0^{-1}$$
is a holomorphic family of $\Psi$DOs of order $\alpha+1$ on $\mathbf{C}\setminus U$, with poles with finite-rank residues at the points of $U$, which coincides with $\Lambda_\alpha^{-1}$ for $\alpha > -1$. Thus
$$H_\alpha(x,y) := \langle Y(\alpha) K_x, K_y \rangle$$
gives, for any $x, y \in \Omega$, a holomorphic function on $\mathbf{C}\setminus U$, with (at most) poles at the points of $U$, which coincides with the weighted harmonic Bergman kernels $H_\alpha(x, y)$ for $\alpha > -1$, proving our claim about the existence of the analytic continuation.

Since the operators $G_\alpha := \mathbf{K} Y(\alpha) \mathbf{K}^*$ still belong to the Boutet de Monvel calculus — being singular Green operators of class zero and order $\alpha$ — the boundary behaviour of the analytically continued kernels $H_\alpha(x, y)$, $\alpha \in \mathbf{C}\setminus U$, is still described by Theorem 10 and Corollary 12.

Finally, let us remark that for $\Omega = \{z \in \mathbf{C} : |z| < 1\}$ the unit disc in $\mathbf{R}^2 \cong \mathbf{C}$ and $w(z) = w(|z|)$ a radial weight, it is elementary that the harmonic Bergman kernels are, up to the constant term, just the real parts of the corresponding holomorphic Bergman kernels:
$$H_w(x,y) = B_w(x,y) + B_w(y,x) - B_w(0,0).$$
Taking in particular $w = w_\alpha$, this implies that the pole-set $U$ above will be the same for the holomorphic and the harmonic kernels. In particular, it follows that the pole-set $U$ can assume the various bizarre forms described in §7.2 of [18].

**7.5 Logarithmic weights.** The results of this paper might extend also to b$\Psi$DOs which are not classical but log-classical (or log-polyhomogeneous), in the sense of allowing symbols with the more complicated expansions

$$k \sim \sum_{j=0}^{\infty} \sum_{m=0}^{M_j} k_{jm}(x', x_n, \xi')$$

where

$$k_{jm}(x', \frac{x_n}{\lambda}, \lambda\xi') = \lambda^{d-j} (\log \lambda)^m k_{jm}(x', x_n, \xi')$$

for $|\xi'| = 1$ and $\lambda \geq 2$; and similarly for the trace and the singular Green symbols. In the language of §7.1, this would lead to the description of e.g. the boundary behaviour of weighted harmonic Bergman kernels $H_w(x, y)$ with weights $w$ of the form

$$w(x) \sim d(x)^\alpha \sum_{j=0}^{\infty} \sum_{m=0}^{M_j} d(x)^j (\log d(x))^m e^{g_{mj}(x)}$$

(where $M_0 = 0$, $\alpha > -1$, $g_{jm} \in C^\infty(\overline{\Omega})$), as has been done for the holomorphic case in [17].

**7.6 Berezin transforms.** The Berezin transform associated with a weighted holomorphic Bergman kernel $B_w(x, y)$ is the integral operator on $\Omega$ defined by

$$\mathcal{B}_w f(x) = B_w(x,x)^{-1} \int_\Omega f(y) |B_w(x,y)|^2 w(y) \, dy.$$

Note that by the reproducing property, $\mathcal{B}_w$ fixes holomorphic and anti-holomorphic functions. For weights $w = w_\alpha$ of the form (82), the asymptotic expansion as $\alpha \to +\infty$ of the Berezin transforms $\mathcal{B}_{w_\alpha} \equiv \mathcal{B}_\alpha$ plays crucial role in quantization of Kähler manifolds (the Berezin and the Berezin-Toeplitz quantizations, see e.g. [1], [38]). It would definitely be of interest to understand the behaviour as $\alpha \to +\infty$ of the similarly defined Berezin transforms also in the harmonic, rather than holomorphic, setting studied in the present paper. While there are a few results available in some special situations, the general case remains unclear.

**7.7 Curvature invariants.** In the holomorphic case, the boundary values of the functions $a, b$ in Fefferman's expansion (1) involve interesting biholomorphic invariants (see e.g. Graham [25]). For the harmonic case, it is clear from the proofs above — and, in fact, is quite immediate from the pseudolocality of b$\Psi$DOs — that the singularity of $K(x, \zeta)$ as some boundary point $x = \zeta = a \in \partial\Omega$, or of $H(x, y)$ at $x = y = a \in \partial\Omega$, is determined completely by the jet of $\partial\Omega$ at $a$. To determine how exactly the Taylor coefficients of the $F(\zeta, r, \nu) \in C^\infty(\partial\Omega \times \overline{\mathbf{R}_+} \times \mathbf{S}^{n-1})$ and $G(x, \zeta) \in C^\infty(\overline{\Omega} \times \partial\Omega)$ in Theorem 13 at $\zeta = a$, $r = 0$ and $x = \zeta = a$, respectively, or of the $F(x, y, r, \nu) \in C^\infty(\overline{\Omega} \times \overline{\Omega} \times \overline{\mathbf{R}_+} \times \mathbf{S}^{n-1})$ and $G(x, y) \in C^\infty(\overline{\Omega} \times \overline{\Omega})$ in Theorem 14 at $x = y = a$, $r = 0$ and $x = y = a$, respectively, depend on this jet — e.g., whether they depend just on some curvature invariants of $\partial\Omega$ at $a$ — seems to be rather difficult. (A hint is given by the formulas (63), (73) in our examples in Sections 5 and 6.) Note that in the holomorphic case, there are many maps that preserve holomorphic functions, i.e. functions annihilated by $\overline{\partial}$ (all biholomorphic mappings); whereas the only maps preserving harmonic functions, i.e. functions annihilated by $\Delta$, in $\mathbf{R}^n$ with $n > 2$ are just the rigid motions and dilations. (See e.g. [3], p. 44.) Consequently, there is no hope of bringing $\partial\Omega$ into some simpler "canonical form" while preserving its harmonic functions (like the Chern-Moser normal form in the holomorphic case). For this reason, it is quite likely that one will need the complete information about the jet of $\partial\Omega$, and not just e.g. the curvature or similar combined quantities, to determine the singularities of the $F$ and $G$ above (or, equivalently, of $K(x,\zeta)$ and $H(x,y)$).

## References

[1] S. Twareque Ali, M. Engliš: *Quantization methods: a guide for physicists and analysts*, Rev. Math. Phys. **17** (2005), 391-490.
[2] S. Axler, P. Bourdon, W. Ramey: *Harmonic function theory,* Springer, 2001.
[3] P. Baird, J.C. Wood, *Harmonic morphisms between Riemannian manifolds,* Clarendon Press, Oxford, 2003.
[4] M. Beals, C. Fefferman, R. Grossman: *Strictly pseudoconvex domains in* $\mathbf{C}^n$, Bull. Amer. Math. Soc. **8** (1983), 125–326.
[5] S. Bergman, M. Schiffer, *Kernel functions and elliptic differential equations in mathematical physics,* Academic Press, New York, 1953.
[6] L. Boutet de Monvel: *Comportement d'un opérateur pseudo-différentiel sur une variété à bord II,* J. d'analyse Math. **17** (1966), 255–304.
[7] L. Boutet de Monvel: *Opérateurs pseudo-différentiels analytiques et problèmes aux limites elliptiques,* Ann. Inst. Fourier (Grenoble) **19** (1969), 169–268.
[8] L. Boutet de Monvel: *Boundary problems for pseudo-differential operators,* Acta Math. **126** (1971), 11–51.
[9] L. Boutet de Monvel, *Pseudo differential operators and their applications,* Duke University, Durham, 1976.
[10] L. Boutet de Monvel, J. Sjöstrand: *Sur la singularité des noyaux de Bergman et de Szegö,* Astérisque **34–35** (1976), 123–164.
[11] A.P. Calderón, *Lecture notes on pseudo-differential operators and elliptic boundary value problems,* Buenos Aires, 1976.
[12] D. Catlin: *The Bergman kernel and a theorem of Tian,* Analysis and geometry in several complex variables (Katata, 1997), Trends in Math., pp. 1–23, Birkhäuser, Boston 1999.
[13] B.R. Choe, H. Koo, H. Yi: *Projections for harmonic Bergman spaces and applications,* J. Funct. Anal. **216** (2004), 388–421.
[14] B.R. Choe, K. Nam: *Berezin transform and Toeplitz operators on harmonic Bergman spaces,* J. Funct. Anal. **257** (2009), 3135–3166.
[15] R.R. Coifman, R. Rochberg: *Representation theorems for Hardy spaces,* Asterisque **77** (1980), 11–66.
[16] R.G. Douglas, *Banach algebra techniques in operator theory,* Academic Press, New York 1972.
[17] M. Engliš: *Weighted Bergman kernels for logarithmic weights,* Pure Appl. Math. Quarterly **6** (2010), 781–813.
[18] M. Engliš: *Analytic continuation of weighted Bergman kernels,* J. Math. Pures Appl. **94** (2010), 622–650.
[19] M. Engliš, D. Lukkassen, J. Peetre, L.-E. Persson: *On the formula of Jacques-Louis Lions for reproducing kernels of harmonic and other functions,* J. reine angew. Math. **570** (2004), 89–129.
[20] C. Fefferman: *The Bergman kernel and biholomorphic mappings of pseudoconvex domains,* Inv. Math. **26** (1974), 1–65.
[21] C. Fefferman: *Parabolic invariant theory in complex analysis,* Adv. Math. **31** (1979), 131–262.
[22] G.B. Folland, *Fourier analysis and its applications,* Wadsworth & Brooks/Cole, Pacific Grove, 1992.
[23] G.B. Folland, *Introduction to partial differential equations. Second edition,* Princeton University Press, Princeton, New Jersey, 1995.
[24] I.C. Gohberg, M.G. Krein, *Introduction to the theory of linear nonselfadjoint operators,* Translations of Mathematical Monographs 18, Amer. Math. Soc., Providence, 1969.

[25] C.R. Graham: *Scalar boundary invariants and the Bergman kernel,* Complex Analysis II (College Park, 1985/86), Lecture Notes in Math. 1276, Springer, Berlin, 1987, pp. 108–135.
[26] G. Grubb, *Distributions and operators,* Springer, 2009.
[27] K. Hirachi: *Invariant theory of the Bergman kernel of strictly pseudoconvex domains,* Sugaku Expositions **17** (2004), 151–169.
[28] L. Hörmander, *The analysis of linear partial differential operators, vol. I,* Grundlehren der mathematischen Wissenschaften, vol. 256, Springer-Verlag, 1985.
[29] L. Hörmander, *The analysis of linear partial differential operators, vol. III,* Grundlehren der mathematischen Wissenschaften, vol. 274, Springer-Verlag, 1985.
[30] M. Jevtić, M. Pavlović: *Harmonic Bergman functions on the unit ball in $\mathbf{R}^n$,* Acta Math. Hungar. **85** (1999), 81–96.
[31] H. Kang, H. Koo: *Estimate of the harmonic Bergman kernel on smooth domains,* J. Funct. Anal. **185** (2001), 220–239.
[32] S.G. Krantz: *Calculation and estimation of the Poisson kernel,* J. Math. Anal. Appl. **302** (2005), 143–148.
[33] J.-L. Lions, E. Magenes, *Problèmes aux limites non homogènes et applications,* vol. 1, Dunod, Paris, 1968.
[34] J. Miao: *Reproducing kernels for harmonic Bergman spaces of the unit ball,* Monatsh. Math. **125**, (1998) 25–35.
[35] K. Nam: *Representations and interpolations of weighted harmonic Bergman functions,* Rocky Mountain J. Math. **36** (2006), 237–263.
[36] U. Neri: *The integrable kernels of certain pseudo-differential operators,* Math. Ann. **186** (1970), 155–162.
[37] S. Rempel, B.-W. Schulze, *Index theory of elliptic boundary problems,* Akademie-Verlag, Berlin, 1982.
[38] M. Schlichenmaier: *Berezin-Toeplitz quantization for compact Kähler manifolds. A review of results,* Adv. Math. Phys. **2010**, Art. ID 927280, 38 pp.
[39] R.T. Seeley: *Complex powers of an elliptic operator,* Singular Integrals, Proc. Symp. Pure Math. X, AMS, Providence, 1967, pp. 288–307.
[40] M.A. Shubin, *Pseudodifferential operators and spectral theory,* Springer-Verlag, Berlin, 2001.
[41] F. Tréves: *Introduction to pseudodifferential and Fourier integral operators,* Plenum, New York, 1980.

Mathematics Institute, Silesian University in Opava, Na Rybníčku 1, 74601 Opava, Czech Republic and Mathematics Institute, Žitná 25, 11567 Prague 1, Czech Republic
*E-mail address*: englis@math.cas.cz